\documentclass
{amsart}

\usepackage{amssymb,mathrsfs}
\usepackage{amsmath}
\usepackage{amsfonts, mathrsfs}
\usepackage{amsthm}
\usepackage{relsize}
\usepackage
{hyperref}

\newcommand{\Be}{\begin{equation}}
\newcommand{\Ee}{\end{equation}}
\newcommand{\Bea}{\begin{eqnarray}}
\newcommand{\Eea}{\end{eqnarray}}
\newcommand{\Bel}{\begin{align}}
\newcommand{\Eel}{\end{align}}
\newcommand{\Beas}{\begin{eqnarray*}}
\newcommand{\Eeas}{\end{eqnarray*}}
\newcommand{\Benu}{\begin{enumerate}}
\newcommand{\Eenu}{\end{enumerate}}
\newcommand{\Bi}{\begin{itemize}}
\newcommand{\Ei}{\end{itemize}}

\theoremstyle{plain}
\newtheorem{thm}{Theorem}[section]
\newtheorem{cor}[thm]{Corollary}
\newtheorem{lem}[thm]{Lemma}
\newtheorem{prop}[thm]{Proposition}

\theoremstyle{remark}
\newtheorem{rmk}{Remark}  
\theoremstyle{definition}
\newtheorem{defn}{Definition}[section]

\numberwithin{equation}{section}
\newcommand{\supp} {\text{\rm supp\! }}
\newcommand{\dist} {\text{\rm dist\! }}

\newcommand{\vep}{\varepsilon}
\newcommand{\classe} {\mathfrak C^D(\varepsilon_\circ)}
\newcommand{\vepc} {\varepsilon_\circ}
\newcommand{\cc}{c_\circ}

\newcommand{\scc}{{s_\circ}}
\newcommand{\os}{\mathcal O_{\!s}}

\begin{document}

\title[maximal estimates for space curves]
{Maximal  estimates for averages 
\\over   space curves 
}

\author[H. Ko]{Hyerim Ko}
\author[S. Lee]{Sanghyuk Lee}
\author[S. Oh]{Sewook  Oh}

\address{Department of Mathematical Sciences and RIM, Seoul National University, Seoul 08826, Republic of Korea}
\email{kohr@snu.ac.kr}
\email{shklee@snu.ac.kr}
\email{dhtpdnr0220@snu.ac.kr}


\begin{abstract}  Let $M$ be the maximal operator associated to  a smooth  curve in 
$\mathbb R^3$ which has nonvanishing curvature and torsion. We prove that $M$  is bounded on $L^p$ if and only if $p>3$.
\end{abstract}

\maketitle 

\section{Introduction}
Let $\gamma$ be a  smooth curve  defined from the interval $\mathbb J:= [-1,1] $ to $\mathbb R^3$. We consider the average $A  f$ over the dilations of $\gamma$ which is given by 
\[
{A  f(x,t)}= \int f(x-t\gamma(s))\psi(s)\,ds, \quad t>0.
\]
Here $\psi$ is a  smooth function with $\supp \psi \subset (-1,1)$. 
We assume that the curve $\gamma$ has  nonvanishing curvature and torsion,  equivalently, 
\begin{equation}\label{nonv}
\det(\gamma'(s), \gamma''(s), \gamma'''(s)) \neq 0
\end{equation}
for $s \in \mathbb J$.  The condition is the natural nondegeneracy condition which is commonly used in  studies related to space curves and the most typical examples are the helix and the moment curve $(s, s^2, s^3)$.   In this paper we are concerned with $L^p$ boundedness of the maximal operator 
\[
M\! f(x)
=\sup_{0<t} \big| A  f(x,t) \big|.
\]

The study of the maximal average over dilated  submanifolds  has  a long history and there is a lot of literature (for example, see \cite{Stein, IKM}  and references therein).   The celebrated Stein's spherical maximal  theorem  \cite{Stein2}  tells that  the spherical maximal function is bounded on $L^p$ if and only if $p>{d}/(d-1)$  for $d\ge 3$.    The case  $d=2$ was later proved by Bourgain \cite{B2}.  As it turned out,  the problem became more difficult  for the circle or the curves with nonvanishing curvature in $\mathbb R^2$ since the typical  interpolation argument relying on  $L^2$ estimate no longer works. 
In such cases the maximal estimates were obtained via the methods of continuum incidence geometry \cite{B2, Sogge, Schlag97, Schlag98} or 
by utilizing  the local smoothing  property  of the averaging operator \cite{MSS, SS, Lee2}.  

Concerning the maximal average over the curve in  three or higher dimensional spaces,  the $L^p$ boundedness  is naturally expected to be even harder to prove since  Fourier transform of the measure supported on  a space curve has slower decay.  $L^p$ boundedness of such maximal operators  has been of interest 
for a long time (see \cite{PS} for a historical comment) but no positive result was known until recently. It was Pramanik and Seeger \cite{PS}  who proved  for the first time that  $M$ is $L^p$ bounded for  $p>38$.  (Also see  \cite{OS, OSS, PS, PS2} for the developments related to $L^p$ Sobolev estimate for the operator $f\to Af(\cdot,t)$.)
Their result was obtained by relying  on Wolff's sharp $\ell^p$ decoupling inequality for the cone in $\mathbb R^3$ \cite{Wolff}. More precisely,  it was shown in \cite{PS} that  the maximal operator $M$ is bounded on $L^p$ for $p>(p_\circ+2)/2$ if 
the sharp $\ell^p$ decoupling inequality holds for $p>p_\circ$. Combined with the  recent $\ell^p$  decoupling inequality on the optimal range $p\ge 6$ which is due to Bourgain and Demeter \cite{BD},  this establishes the $L^p$ boundedness for $p>4$.  However,  a modification  of  Stein's example in  \cite{Stein2} shows  that   $M $ can not be bounded on $L^p$ for $p\le 3$ (see Section  \ref{Sec:N}).

In this paper we fill the gap and settle the problem of $L^p$ boundedness of $M$.

\begin{thm}\label{max}
Suppose that $\gamma :\mathbb J \rightarrow \mathbb R^3$ is a smooth curve which has nonvanishing curvature and torsion, and  $\psi$ is a nontrivial, nonnegative, smooth function supported in $(-1,1)$. Then, there is a constant $C$ such that  
\begin{align}\label{max-est}
\| M\! f\|_{L^p(\mathbb R^3)}
\le C \|f\|_{L^p(\mathbb R^3)}
\end{align}
for all $f\in L^p(\mathbb R^3)$ if and only if $p>3$.
\end{thm}
  
The assumption that $\psi$ is smooth is not necessary and it is clear that   the theorem   holds true for a continuous $\psi$.  Even though $\gamma$ is assumed to be smooth, there is a positive integer $D$ such that    \eqref{max-est} holds for  $\gamma \in \mathrm C^D(\mathbb J)$ (see Remark \ref{order-} at the end of  Section \ref{sec3}).  

The maximal estimate  in \cite{PS} was shown by exploiting  $L^p$ local smoothing phenomena  of the averaging operator. However, compared with  the average over hypersurfaces or curves in $\mathbb R^2$, the $L^p$ local smoothing property of $A$ is not  well understood.  
We instead try to
make use of an $L^p$-$L^q$ type smoothing estimate 
which has a close connection to the adjoint restriction estimate. Usefulness of such estimates has been manifested  in the study of $L^p$ improving property of the localized circular and spherical maximal functions \cite{SS, Lee2} (also see \cite{AHRS, RS, BORSS}).   

Our argument in this paper  is closely related  to  the induction strategy developed 
by Ham and one of the authors  \cite{HL}. They  obtained  the sharp  adjoint restriction estimate for the space curve in $L^p(\mu)$ when $\mu$ is an $\alpha$-dimensional measure (see  Section \ref{sec2.1} for the definition). The work   was in turn inspired by the multilinear approach due to Bourgain and Guth \cite{BG}.  
Main novelty of  the current paper lies in devising an induction argument which directly works for the  maximal  operator.  
In contrast to the adjoint restriction operator a suitable form of multilinear estimate is not so obvious for the averaging operator $A$.  
In order to prove  a multilinear  
estimate for $A$ which enjoys  a better boundedness property under a certain additional assumption,  we first express the operator $A$ as a sum of adjoint restriction operators and then relate them to geometry of the curves so that  the transversality condition 
can be reformulated  in terms of the relative positions between the associated curves. 
Unfortunately, some the consequent adjoint restriction operators are associated to $\mathrm C^{1,1/2}$ surfaces  but  not to $\mathrm C^2$  surfaces, so we can not directly apply  the multilinear restriction estimate which is due to  Bennett, Carbery, and Tao \cite{BCT}. 
However, 
it is not difficult to see that the argument in \cite{BCT}  continues to work for the $\mathrm C^{1,1/2}$ surfaces  (see Theorem \ref{holder-rest} below).  We also make use of some of the results from \cite{PS} to strengthen the multilinear estimate and also to deal with the nondegenerate part, whereas the difficult degenerate part is to be handled by  the multilinear estimate which we prove in  Section \ref{sec3}. 

The argument here can be further developed  to prove not only $L^p$ improving property of the maximal operator $ \sup_{1\le t\le 2} \big| A  f(x,t) \big|$ but also maximal estimates with respect to $\alpha$-dimensional measures (see Remark \ref{further-result} at the end of Section \ref{s-conclusion}). 
Nonetheless,  we do not attempt to pursue the matter in this paper. 

\medskip

\noindent{\emph{Structure of the paper.}}
In Section \ref{sec2} we show that  the maximal estimate can be  deduced from a form of weighted estimate, and 
we formalize the induction setup to prove the weighted estimate. 
 In Section \ref{sec3} we obtain a weighted multilinear  estimate for $A$ under a certain  separation condition. In Section \ref{sec4} we  establish the maximal bound putting  
 the previous estimates together and show the optimality of the range  of $p$.

\section{Reductions and preliminaries }\label{sec2}
In this section 
we reduce the proof of maximal estimate to showing  a form of weighted estimates for the averaging operators which are given   by the curves close to a specific curve.  We also obtain some preparatory results which are to be used to prove  the estimates in Section 3 and Section 4.

By the argument  in \cite{B2} (also see  \cite{Schlag96}), which relies on  Littlewood-Paley decomposition and scaling,   one can obtain
 the maximal estimate \eqref{max-est} from 
that for $\sup_{1\le t\le 2} \big| A  f(x,t) \big|$.  More precisely,  it is sufficient  to show that  there is an  $\varepsilon_p>0$ such that  
\begin{equation}\label{Lmax}
\| A  f\|_{L_x^pL_t^\infty(\mathbb R^3 \times [1,2])}
\le C \lambda^{-\varepsilon_p} \| f\|_{L^p(\mathbb R^3)}
\end{equation}
holds for  all $f\in \mathcal S(\mathbb R^3)$ whenever   
\Be \label{lambda-f} \supp \widehat f \subset \mathbb A_\lambda:=\{ \xi \in \mathbb R^3:
3 \lambda/4 \le |\xi| \le 7\lambda/4 \}, \ \ \lambda \ge1.\Ee
For the rest of the paper, we assume \eqref{lambda-f} unless it is mentioned otherwise.  

\medskip

\noindent{\bf Notation.}  
Throughout the paper   $C$, $C_1,\dots$ and $c$ are supposed to be  independent positive constants, and $C_\varepsilon$, $C_\delta$ are constants depending  on $\varepsilon, \delta$ but all of these constants may vary at each appearance.  In addition to the conventional notation $\widehat \cdot$\, we  use $\mathcal F$ and $\mathcal F^{-1}$ to denote the Fourier and inverse Fourier transforms, respectively.  By $Q_1=\mathcal O(Q_2)$ we denote $|Q_1|\le  CQ_2$ for a constant  $C$ and  we also use the notation  $Q_1=\mathcal O_{\!s}(Q_2)$ if   $|Q_1|\le Q_2$. 

\subsection{Estimate with $\alpha$-dimensional measure} 
\label{sec2.1} Let $\mathbb B^d(z,r)$ denote the ball of radius $r$ which is centered at $z\in \mathbb R^d$. Let $\mu$  be a positive Borel measure on $\mathbb R^4$. 
For $0<\alpha \le 4$  we say  $\mu$ is $\alpha$-dimensional if there is a constant $C$ such that
\[
\mu (\mathbb B^4(z,r))
 \le C r^\alpha
 \]
for all $r>0$ and $z \in \mathbb R^4$.  
For an $\alpha$-dimensional measure $\mu$ we define 
\[\langle \mu\rangle_\alpha = 
\sup_{z \in \mathbb R^4,\, r>0} r^{-\alpha} \mu (\mathbb B^4(z,r)).\]
 Instead of directly proving the maximal estimate  \eqref{Lmax} we obtain estimates for $Af$  with  $\alpha$-dimensional measures.  From those estimates  we can deduce the  estimate \eqref{Lmax}. As far as the authors are aware, 
 it seems that 
 this type of argument deducing the maximal estimate  from the estimates with  $\alpha$-dimensional measures 
   first  appeared in \cite{OO}. (See  also \cite[p.1283]{Wolff} for a related discussion.)

\begin{thm}\label{LSA} 
 Let $\mu$ be 3-dimensional. Suppose that
$\gamma: \mathbb J \rightarrow \mathbb R^3$ is a smooth curve satisfying \eqref{nonv}. Then, for  $p>3$ there is an $\varepsilon_p>0$ such that
\begin{align}\label{smoothing}
\| A   f\|_{L^p(\mathbb R^3 \times [1,2],d\mu)}
\le C \langle \mu \rangle_3^{\frac 1p} \lambda^{-\varepsilon_p} \| f\|_{L^p(\mathbb R^3)}
\end{align}
holds whenever $\widehat f$ is supported on $\mathbb A_\lambda$. 
\end{thm}

We shall work only with $3$-dimensional  measures even though it is possible to prove such  estimates  with $\alpha$-dimensional measure, $\alpha\neq 3$ on a certain range of $p$ (see Remark \ref{further-result}).  The following  shows
the estimate \eqref{smoothing} implies \eqref{Lmax}.

\begin{lem}
\label{measure-max}
Suppose \eqref{smoothing} holds true for all $3$-dimensional measures $\mu$. Then the estimate \eqref{Lmax} holds. 
\end{lem}

\noindent  To prove  this, we start with  an elementary lemma. 
\begin{lem} \label{lem:ker} Let  $\eta \in \mathrm C_0^\infty([2^{-3},  2^3])$ and $\psi\in \mathrm C_0^\infty(\mathbb J)$.  Set  $r_0=  1+ 4\max\{|\gamma(s)|: s \in {\rm\supp} \psi\} $ and
  \[
 K_\eta (x,t)=(2\pi)^{-3} \iint e^{i(x\cdot \xi -t\gamma(s)\cdot \xi)}\psi(s)\,ds \, \eta(\lambda^{-1}|\xi|)\,d\xi.
\]
If $|x|\ge r_0$ and $|t|\le 2$, then $ | K_\eta(x,t)| \le    C   \|\eta\|_{\mathrm C^{2N+3}} E_N(x)
$ for any $N\ge1$   where $E_N(x):=\lambda^{-N}(1+|x|)^{-N}$. 
\end{lem}

\begin{proof}
 We see 
 $
K_\eta (x,t)=\frac{\lambda^3}{(2\pi)^3}\iint e^{i\lambda(x\cdot \xi -t\gamma(s)\cdot \xi)}\psi(s)\,ds \, \eta(|\xi|)\,d\xi 
$
by  changing variables $\xi\to \lambda\xi$. Then repeated integration by parts in $\xi$  gives the desired estimate since $|\nabla_\xi(x\cdot \xi-t\gamma(s)\cdot \xi)|\ge  {2^{-1}|x|}$ if $|x|\ge r_0$ and $|t|\le 2$.
\end{proof}

\begin{proof}[Proof of Lemma \ref{measure-max}]
To obtain  \eqref{Lmax}  it suffices to show the local estimate 
\begin{align}\label{Lmax**}
\| A  f\|_{L_x^pL_t^\infty(\mathbb B^3(0,1) \times [1,2])}
\le C \lambda^{-\varepsilon_p} \| f\|_{L^p(\mathbb R^3)}.
\end{align}
This is obvious if   $\widehat f$ is not assumed to be  supported in $\mathbb A_\lambda$.  
However,  we may  handle $f$ as if it were supported on a ball of radius $r_0$. 
Since $\supp \widehat f\subset \mathbb A_\lambda$, $A  f(\cdot,t)=K_\eta (\cdot,t)\ast f$ for an $\eta$ such that  
 $\eta \in \mathrm C_c^\infty((2^{-1},  2))$ and   $\eta=1$  on $[3/4, 7/4]$. So, 
 Lemma \ref{lem:ker}  gives $|K_\eta (x,t)| \le C E_N(x)$ if $|x|\ge r_0$ and $|t|\le 2$. Thus, by the  typical localization argument 
 (e.g., see the proof of Lemma \ref{localtoglobal})  one can easily see that \eqref{Lmax**} implies \eqref{Lmax}.

In order to prove  \eqref{Lmax**},  using the Kolmogorov-Seliverstov-Plessner linearization,   it is enough   to show 
\begin{equation}\label{msr1}
\|A   f(\cdot, {\bf{t}}(\cdot)) \|_{L^p(\mathbb B^3(0,1))} \le C
\lambda^{-\varepsilon_p} \|f\|_{L^p(\mathbb R^3)}
\end{equation} 
for all measurable function ${\bf{t}}: \mathbb B^3(0,1)\to [1,2]$ with $C$ independent of ${\bf{t}}$. Since $\widehat f$ is supported in $\mathbb A_\lambda$, $A  f$ is uniformly continuous on every compact subset.  So,  for  \eqref{Lmax**} it is sufficient to show \eqref{msr1} while  assuming ${\bf{t}}$ is continuous. With a continuous function 
${\bf{t}}$, the positive linear functional  $\mathrm C_c(\mathbb R^4)\ni F \mapsto \int_{\mathbb B^3(0,1)} F(x,{\bf{t}}(x))dx$  defines  a  measure $\mu$\footnote{In fact,  we see $\mu$ is  a regular Borel measure by the Riesz-Markov-Kakutani representation theorem.}  by the relation
\begin{align*}
\int F(x,t)\, d\mu(x,t)
=\int_{\mathbb B^3(0,1)} F(x,{\bf{t}}(x))\,dx, \quad  F \in \mathrm C_c(\mathbb R^4). 
\end{align*}
We now notice  that $\mu$ is a $3$-dimensional measure. 
 Since  
 $\mathbb B^4((x_\circ,t_\circ),r) \subset \{ (x,t)  \in \mathbb R^3\times \mathbb R: |x-x_\circ| \le r\}$, 
\begin{align*}
\mu \big( \mathbb B^4((x_\circ,t_\circ),r) \big)
=
\int_{\mathbb B^3(0,1)} \chi_{\mathbb B^4((x_\circ,t_\circ),r)}(x,{\bf{t}}(x))\,dx
\le
\int \chi_{\mathbb B^3(x_\circ,r)}(x)\,dx = \frac 43 \pi r^3
\end{align*}
for any $r>0$ and $(x_\circ, t_\circ)\in \mathbb R^3\times \mathbb R$. Thus we have $\langle \mu \rangle_3\le   4 \pi/3.$ 
Noting  $\|A  f(\cdot, {\bf{t}}(\cdot)) \|_{L^p(\mathbb B^3(0,1))}=\|A  f \|_{L^p(d\mu)}$,  
we apply Theorem \ref{LSA} and get  \eqref{msr1} with $C$ independent of  ${\bf{t}}$. 
\end{proof}

\subsection{Weighted estimate} 
For $0<\alpha \le 4$,   let us denote  by $\Omega^\alpha$ the collection of nonnegative  measurable functions $\omega$ on $\mathbb R^4$ such that 
the measure $ \omega\, dxdt$ is $\alpha$-dimensional. For a simpler notation  we 
denote 
\[
[\omega ]_\alpha 
=\langle  \omega\,dxdt  \rangle_\alpha.
\]  
Even though $\Omega^\alpha$ is properly contained in the set of  $\alpha$-dimensional measures,  
the fact that $\supp \widehat f\subset \mathbb A_\lambda$  allows us to recover  
the  estimate \eqref{smoothing}  from an estimate against $\omega \in \Omega^\alpha$.  

\begin{lem}\label{WtoM}
Let $I=[2^{-1},2^2]$. 
Suppose that  
\begin{align}\label{Weighted}
\| A   f\|_{L^p(\mathbb R^3 \times I, \omega)}
\le C [ \omega ]_3^{\frac 1p} \lambda^{-\varepsilon_p} \|f\|_{L^p(\mathbb R^3)}
\end{align}
holds whenever $\omega \in \Omega^3$ and $\widehat f$ is supported on $\mathbb A_\lambda$. Then \eqref{smoothing} holds for any  $3$-dimensional measure $\mu$.
\end{lem}

The proof of the maximal estimate \eqref{Lmax} is now  reduced to showing \eqref{Weighted}.  Lemma \ref{WtoM}  of course remains valid for any $\alpha\in (0,4]$.  

To show  Lemma \ref{WtoM} we make use of the next two lemmas: Lemma \ref{omega-omega} and \ref{lem:ker2}. The former 
can be shown  following  the standard argument (for example, see \cite[pp. 47--49]{Matt}), so we omit the proof.

\begin{lem}
\label{omega-omega}
Let $0<\alpha\le 4$ and  $\varphi \in \mathcal S(\mathbb R^4)$. Set
$\varphi_\lambda=\lambda^4\varphi(\lambda\, \cdot)$. 
If $\mu$ is an $\alpha$-dimensional measure,  then $|\varphi|_\lambda \ast \mu\in \Omega^\alpha$ and $[|\varphi|_\lambda \ast \mu]_\alpha \le C_\varphi \langle \mu \rangle_\alpha$. 
\end{lem}

\noindent

In what follows,  $\widetilde \chi$ denotes a function in $\mathrm C_0^\infty(I) $ which satisfies
$\widetilde \chi =1$ on $[1,2]$,  and $\beta$,  $\beta_0$ respectively denote the  functions such that $\beta \in \mathrm C_0^\infty([2^{-1},  2])$, $\beta=1$  on $[3/4, 7/4]$; $\beta_0 \in {\mathrm C_0^{\infty}} ([-2,2])$,  $\beta_0=1$ on $[-1,1]$.

\begin{lem}\label{lem:ker2}  Let $r_0=  1+ 4\max\{|\gamma(s)|: s \in {\rm\supp} \psi\} $  and let
\[
m (\xi, \tau) =   \iint \widetilde \chi(t)\,  e^{-i t(\tau+ \gamma(s)\cdot \xi)} \psi(s)\,ds dt\, {\beta}(\lambda^{-1} |\xi|), \quad (\xi, \tau)\in \mathbb R^3\times \mathbb R. 
\]
Then, we have 
$
|\mathcal F^{-1}\big(m(\xi, \tau)( 1- \beta_0 ((\lambda r_0)^{-1} \tau))\big)|\le C_N\| \psi\|_{\infty} \widetilde E_t^N$  for any $N>0$  where $\widetilde E_t^N:=(1+|t|)^{-N} E_N$.  
\end{lem}

\begin{proof} Let $\rho_\ell (t)=(-it)^{k+\ell}\widetilde\chi(t)$ and note  that $ \partial_\xi^\alpha  \partial_\tau^{k} m(\xi, \tau)$ is  a sum  of the intergrals  
$
\int \widehat \rho_{|\alpha_1|} (\tau+ \gamma(s)\cdot \xi) (\gamma (s))^{\alpha_1}\psi(s)\,ds
\times \mathcal O(\lambda^{-|\alpha_2|})
$ 
with $\alpha_1+\alpha_2=\alpha$. Thus it follows that $| \partial_\xi^\alpha  \partial_\tau^{k} m(\xi, \tau) |\le C_N \| \psi\|_{\infty}  r_0^{|\alpha|}  (r_0\lambda)^{-N}(1+|\tau|)^{-N}$ for any $N$ if $|\tau|\ge r_0\lambda$. We then get the desired estimate  by routine integration by parts.
 \end{proof}

\begin{proof}[Proof of Lemma \ref{WtoM}]
We define an auxiliary operator 
$\widetilde  A $ 
by 
\[
\mathcal F(\widetilde  A h)(\xi,\tau)
=
\beta_0 ((\lambda r_0)^{-1} \tau)\mathcal F\big(\widetilde \chi(t) A h\big)(\xi,\tau).
\]
 Since $\widehat f$ is supported in $\mathbb A_\lambda$,   we have
 $  |(\widetilde \chi(t) A  - \widetilde  A) f|\le C\widetilde E_t^N \ast |f|$ by Lemma  \ref{lem:ker2}. 
 We then note  that $\int \widetilde E_t^N(x-y) d\mu(x,t) \le C  \lambda^{-N} \langle \mu \rangle_3$ and  $\int \widetilde E_t^N(x-y) dy \le C \lambda^{-N}$. 
 Thus by Schur's test  we get
\begin{align}\label{At-A}
\| \widetilde E^N_t \ast h
\|_{L^p(\mathbb R^3 \times \mathbb R, d\mu)} 
\le C \langle \mu \rangle_3^{\frac 1p} \lambda^{-N} \|h\|_{L^p(\mathbb R^3)}
\end{align}
  for $1 \le p \le \infty$ and a large $N$. So, in order to show \eqref{smoothing},
  it suffices   to prove
\begin{align}\label{claim1}
\| \widetilde  A f\|_{L^p(\mathbb R^3 \times [1,2], d\mu)}
\le C \langle \mu \rangle_3^{\frac 1p} \lambda^{-\varepsilon_p} \|f\|_{L^p(\mathbb R^3)}.
\end{align}

Since  the space time Fourier transform of $\widetilde  A f$ is supported  in $\mathbb B^4(0,2^2r_0\lambda)$, 
${\widetilde  A f}={\widetilde  A f}\ast \varphi_{r_0\lambda}$ for some $\varphi\in\mathcal S(\mathbb R^4)$. This gives $|{\widetilde  A f}|^p\le C  |\widetilde  A f|^p  \ast |\varphi_{r_0\lambda}|$ via  H\"older's inequality.  Thus  we have 
\[
\|\widetilde  A f\|_{L^p(\mathbb R^3 \times [1,2], d\mu)}
\le C 
\| \widetilde  A f\|_{L^p(\mathbb R^3 \times \mathbb R, \omega)},
\]
where $\omega=|\varphi_{r_0\lambda}|\ast \mu$. Therefore, using $ |(\widetilde \chi(t) A  - \widetilde  A) f|\le C\widetilde E_t^N \ast |f|$ again, 
we  have only to obtain the estimate for $\widetilde \chi(t) A  f$ 
in $L^p(\mathbb R^3 \times \mathbb R, \omega)$ since the minor part can be handled as before. Since $[\omega]_3\le C\langle\mu\rangle_3$ by Lemma \ref{omega-omega}, 
the estimate  \eqref{claim1} follows from \eqref{Weighted} because $\supp \widetilde \chi\subset  I$.
\end{proof}

\subsection{Normalization of curves and weights}
In order to prove the  estimate \eqref{Weighted},  as mentioned before,   we use an induction type argument over a class of curves. For the purpose  we need to normalize the curves properly so that the induction assumption applies. 
This step is especially  important for defining  the induction quantity and proving uniform estimates 
(cf. \cite{HL, Lee}).

Let $D\ge 2^5$ be a positive integer which is taken to be large. Let $\gamma \in \mathrm C^D(\mathbb J)$ which satisfies \eqref{nonv}.  Then, 
for $ s_\circ $ and $0<\delta \ll1$ such that $[s_\circ-\delta,s_\circ+\delta] \subset \mathbb J$,  we  define 
\[ 
\mathrm M_{\gamma}^\delta (s_\circ)=  
\big( \delta \gamma'(s_\circ), \delta^2 \gamma''(s_\circ), \delta^3 \gamma'''(s_\circ)
\big)
\]
and 
\begin{align}
\label{scaledcurve}
\gamma_{s_\circ}^\delta(s) =
(\mathrm M_{\gamma}^\delta (s_\circ))^{-1}
\big( \gamma(\delta s+s_\circ)-\gamma(s_\circ) \big).
\end{align}

Let  $\gamma_\circ(s)=(s, s^2/2!, s^3/3!)$.  We  consider a class of curves  
which are small perturbations of the curve $\gamma_\circ$ in $\mathrm C^D(\mathbb J)$.  
For $\varepsilon_\circ>0$,   we set
\[
\mathfrak C^D(\varepsilon_\circ)
=\big\{ \gamma \in \mathrm C^D(\mathbb J) :  \ \| \gamma- \gamma_\circ \|_{\mathrm C^D(\mathbb J)} \le \varepsilon_\circ \big \}.
\]
Using an affine map,  one can transform a small enough sub-curve of any $\gamma\in \mathrm C^D(\mathbb J)$ satisfying \eqref{nonv}  so as  to be contained in $\mathfrak C^D(\varepsilon_\circ)$.   The following lemma is a slight modification of  \cite[Lemma 2.1]{HL}. 

\begin{lem}\label{lem:stable} 
Let  $s_\circ \in (-1,1)$ and $\gamma\in \mathrm C^D(\mathbb J)$ satisfy \eqref{nonv} on $\mathbb J$. Then, 
for  any $\varepsilon_\circ>0$, there exists $\delta_\ast=\delta_\ast(\varepsilon_\circ,\gamma)>0$ such that
$\gamma_{s_\circ}^\delta \in \mathfrak C^D(\varepsilon_\circ)$ whenever
 $[s_\circ-\delta,s_\circ+\delta] \subset \mathbb J$ and $|\delta| \le \delta_\ast$. Additionally, if $\gamma\in \mathfrak C^D(\varepsilon_\circ)$ and $\vepc<2^{-5}$,  then
  there is  a  uniform  $\delta_\circ>0$ such that 
$\gamma_{s_\circ}^\delta \in \mathfrak C^D(\varepsilon_\circ)$ whenever  $[s_\circ-\delta,s_\circ+\delta] \subset \mathbb J$ with $|\delta| \le \delta_\circ$.
\end{lem}

For a matrix $\mathrm M$ we denote $\|\mathrm M\|=\sup_{|z|=1} |\mathrm Mz|$.

\begin{proof} 
By Taylor expansion of $\gamma(\delta \cdot +s_\circ)$ about $s=0$,
we have
\begin{align*}
\gamma(\delta s+s_\circ)-\gamma(s_\circ)
&= \delta \gamma'(s_\circ) s + \delta^2 \gamma''(s_\circ) \frac{s^2}{2} 
+\delta^3 \gamma'''(s_\circ) \frac {s^3} {3!}+ \widetilde R(s_\circ,\delta,s)\\
&=\mathrm M_{\gamma}^\delta (s_\circ) \gamma_\circ(s)+\widetilde R(s_\circ,\delta,s)
\end{align*}
and $\| \widetilde R(s_\circ,\delta,\cdot )\|_{\mathrm C^D(\mathbb J)} \le C\delta^4$.
By \eqref{scaledcurve},  $
\gamma_{s_\circ}^\delta(s)=\gamma_\circ(s)+ (\mathrm M_{\gamma}^\delta (s_\circ))^{-1} \widetilde R(s_\circ,\delta,s).
$
Since $\|(\mathrm M_{\gamma}^\delta (s_\circ))^{-1} \|\le C_1\delta^{-3}$ for a constant $C_1$, 
taking  a positive  $\delta_\ast$  such that  $CC_1 \delta_\ast\le \vepc$ we have 
$\|(\mathrm M_{\gamma}^\delta (s_\circ))^{-1} \widetilde R(s_\circ,\delta,\cdot)\|_{\mathrm C^D(\mathbb J)}
 \le \varepsilon_\circ$ and, hence, $\gamma_{s_\circ}^\delta \in \mathfrak C^D(\varepsilon_\circ)$ for $0<\delta\le\delta_\ast$.  The second assertion can also be shown 
in the same manner, so we omit the detail. 
\end{proof}

For $\delta>0 $   we 
denote  by  $\mathrm D_\delta $ the diagonal  matrix  
$(\delta e_1,\delta^2 e_2,\delta^3 e_3).$ 
To normalize the weights we need  the next lemma, which one can show  by following  the argument in  \cite{HL}.

\begin{lem}
\label{weight}
Let $0< \alpha \le 4$,  $0<\delta \ll1$ and $\omega\in \Omega^\alpha$, and  let $\mathrm M$ be a  $4 \times 4$ nonsingular matrix. Set
$
\omega^\delta(x,t)
=\omega \big( \mathrm D_\delta x,t \big)
$
and $\omega_{\mathrm M}(x,t)
=\omega \big( \mathrm M( x,t) \big). $
Then, for a constant $C$ independent of $\omega$ and $\delta$,   we have
\begin{align}
\label{scaling-weight}
[\omega^\delta]_\alpha
&\le C \delta^{3\alpha-12} [\omega]_\alpha, 
\\
\label{scaling-weight2}
[\omega_{\mathrm M}]_\alpha
&\le |\!\det\mathrm M|^{-1} \|\mathrm M\|^\alpha   [\omega]_\alpha.
\end{align} 
\end{lem}

\begin{proof} The inequality \eqref{scaling-weight} is equivalent to  
\[
\int_{\mathbb B^4(y,r)} \omega(\mathrm D_\delta x,t)\,dxdt
\le C \delta^{3\alpha-12} [\omega]_\alpha   r^\alpha
\]
 for $y\in\mathbb R^4$ and $r>0$.  Changing  variables $x\to  \mathrm D_\delta^{-1}x$,  we see  the left hand side is equal to  $\delta^{-6}\int \chi_{\mathbb B^4(y,r)}(\mathrm D_\delta^{-1}x,t)\,\omega(x,t) \,dxdt$. Then we note  that  the set   $\{(x,t): (\mathrm D_\delta^{-1} x, t) \in \mathbb B^4(y,r)\}$ is contained in a rectangle $\mathcal R_\delta$ of dimensions  about
$\delta r \times \delta^2 r \times \delta^3 r \times r$.
Since $\mathcal R_\delta$ is  covered by
at most $C\delta^{-6}$ many  balls  of radius $\delta^3r$,  \eqref{scaling-weight} follows.   

For \eqref{scaling-weight2} we only have to show   
\[
\int_{\mathbb B^4(y,r)} \omega(\mathrm M(x,t))\,dxdt
\le  |\!\det\mathrm M|^{-1}\| \mathrm M\|^{\alpha}  [\omega]_\alpha  r^\alpha
\]
 for $y\in \mathbb R^4$ and $r>0$.  Changing variables, we see that  the left hand side  equals $  |\!\det\mathrm M|^{-1}\int \chi_{\mathbb B^4(y,r)} (\mathrm M^{-1}(x,t))\omega(x,t)  dxdt$. So, we get  the inequality \eqref{scaling-weight2} since  $(x,t) \in  \mathbb B^4(\mathrm My,\|\mathrm M\|r)$ if $\mathrm M^{-1}(x,t) \in \mathbb B^4(y,r)$. 
\end{proof}

\subsection{Reduction and the induction quantity}  
\label{sec2.4}
Throughout  the paper we fix a small  positive constant  $c_\circ$. To show \eqref{Weighted} for a smooth curve satisfying \eqref{nonv}, it is sufficient to handle  
$\gamma \in \mathfrak C^D(\varepsilon_\circ)$ with a small $\varepsilon_\circ>0$ while $\psi \in \mathrm C^D$ and  $\supp \psi \subset 
[-c_\circ,c_\circ]$.  As we shall see later,  this can be shown 
by a finite decomposition and changing variables via affine transformations.

\begin{defn} 
Let $\cc$, $\vepc$ and $\delta$ be the numbers such that    $ 0<c_\circ\le 2^{-10}$,  $ 0<\vepc\le c_\circ^2$, and 
\Be  \label{del-choice}
0<\delta\le \min(\cc, \delta_\circ)
\Ee
 where 
$\delta_\circ$ is given in Lemma \ref{lem:stable}. The number $\delta$ is  to be 
chosen later (see Section \ref{s-conclusion}). We  denote $J_\circ=[-c_\circ,c_\circ]$, and we set 
\[ \mathfrak J(\delta)=\big\{ J: J=[c_\circ\delta (k-1), c_\circ\delta (k+1)], \, k\in \mathbb Z, \,  |k|\le    (\cc\delta)^{-1} +1\big\} \] 
so that  the intervals in $\mathfrak J(\delta)$ cover $\mathbb J$.  For each $J\in \mathfrak J(\delta)$ we define   $\mathfrak N^D(J)$ to be  the set of functions such that 
 $\psi \in \mathrm C^D_0(J)$ and $\|  \psi(|J|\,\cdot) \|_{\mathrm C^D(\mathbb R)}\le  1$. For a given interval  $J$  we denote by $\psi_{\!J}$ a function in  $\mathfrak N^D(J)$. 
 \end{defn}

For a smooth function $a$ on $\mathbb J \times I \times
\mathbb A_\lambda$, following \cite{PS}, we  define an integral operator by setting 
\begin{align}\label{Aaf}
A^\gamma  [a]f(x,t)
=(2\pi)^{-3}\iint e^{i (x -t\gamma(s))\cdot \xi} 
a(s,t,\xi)\,ds \widehat f(\xi)\,d\xi.
\end{align}
In particular, we note  $A  f=A^\gamma [\psi]f$ as is clear by Fourier inversion.

Let us take $\zeta\in \mathrm C_0^\infty([-1, 1])$ such that   $\zeta \ge 0$ and $\sum_{k\in\mathbb Z} \zeta(s-k)=1$. For an interval $J$ we denote by $c_J$ the center of  $J$ and set  
$\zeta_J(s)= \zeta(2(s-c_J)/|J|)$. Consequentially, $\zeta_J\in \mathrm C_0^\infty(J)$ and 
$
\sum_{J\in \mathfrak J(\delta)} \zeta_J(s)=1$ for $s\in \mathbb J. 
$
As a result, we have  
\begin{align}\label{Jsum}
A^\gamma [\psi] f(x,t) =\sum_{J\in \mathfrak J(\delta)} A^\gamma [\psi \zeta_J]  f(x,t)
\end{align}
if $\supp \psi\subset \mathbb J$. 
The following is one of the key lemmas which relates the estimate for the average over a short curve to that  over a larger one. 
\begin{lem}\label{lem:scaled}
Let $I'\subset I $ be an interval, and let $\omega\in \Omega^3$, $J=[\scc-\cc\delta, \scc+\cc\delta]\in \mathfrak J(\delta)$   
and  $\psi_{\!J}\in \mathfrak N^D(J)$. 
Suppose  that $\gamma\in \mathrm C^D(\mathbb J)$  satisfies  \eqref{nonv} 
and  $\supp \widehat  f\subset \mathbb A_\lambda$. Then,  there are  $\widetilde \omega\in \Omega^3$,  $\widetilde f$ with $\|\widetilde f\|_p=\|f\|_p$, and  $\psi_{\!J_\circ} \in\mathfrak N^D(J_\circ)$ which satisfy the following: 
\begin{align}
\label{scale3} 
&\| A^\gamma [\psi_{\!J}] f\|_{L^p(\mathbb R^3 \times I', \omega)} = 
\delta^{1-\frac{3}p}  \| A^{ \gamma_{\scc}^\delta} [ \psi_{\!J_\circ}] \widetilde f\|_{L^p(\mathbb R^3 \times I', \widetilde \omega)},
\\
\label{weight-}
&[ \widetilde \omega]_3\le  C(1+|\gamma(s_\circ)|)^3  |\det \mathrm M_{\gamma} ^1(s_\circ)|^{-1} (1+\|\mathrm M_{\gamma}^1 (s_\circ)\|)^3 [\omega]_3,
\end{align}
and
\begin{eqnarray}
\label{supp-}
\hspace{42pt}\supp \mathcal F(\widetilde  f\,)\subset  \Big\{\xi: \frac34 d_\ast \delta^3\lambda\le |\xi|\le \frac 74 d^\ast \delta\lambda\Big\},
\end{eqnarray} 
where $1/d_\ast= \|(\mathrm M_{\gamma}^1 (s_\circ))^{-t}\|$ and $1/d^\ast =\inf_{|z|=1} |(\mathrm M_{\gamma}^1 (s_\circ))^{-t} z|$.  
\end{lem}

\begin{proof}
We denote  $\psi_{\!J_\circ}(s)=\psi_{\!J}(\delta s+\scc)$. It is clear that $\psi_{\!J_\circ} \!\!\in  \mathfrak N^D(J_\circ)$.  We  set 
\[  \widetilde f(x)=  |\!\det ( \mathrm M_{\gamma}^\delta(\scc))|^{\frac1p} 
f( \mathrm M_{\gamma}^\delta(\scc) x).\] 
Then,  we see   $\|\widetilde f\|_p=\|f\|_p$ and  $\mathcal F(\widetilde  f\,)$ is supported in the set  
 $S_\lambda=\{\xi:  3\lambda/4\le  |(\mathrm M_{\gamma}^\delta(\scc))^{-t}\xi|    \le 7 \lambda/4 \}$ because 
 $\supp \widehat f\subset  \mathbb A_\lambda$.  
 Since $\mathrm M_{\gamma}^\delta (s_\circ)=\mathrm M_{\gamma}^1(\scc) \mathrm D_\delta$, 
 it is easy to see that 
$S_\lambda\subset \{\xi:  3\lambda d_\ast /4\le  |\mathrm D_\delta^{-1} \xi|    \le 7 \lambda d^\ast/4 \}$, thus we get \eqref{supp-}.

 We now define  $\overline \omega$ 
and $\widetilde \omega$  by setting $\overline \omega(x,t)= \omega (x+t\gamma(\scc),t)$ and 
\[ 
\widetilde \omega(x,t)=\delta^{3}\overline\omega(\mathrm M_{\gamma}^\delta(\scc) x,t),
\]
respectively.  Denoting  by $\mathrm M$  the matrix such that  $\mathrm M(x,t)=(x+t\gamma(\scc),t)$,   we note that 
$\overline \omega=\omega_{\mathrm M}$, 
$\det \mathrm M=1$, and $\|\mathrm M\|\le 1+|\gamma(s_\circ)|$. Thus using  \eqref{scaling-weight2} 
we  have
$[\overline \omega]_3\le (1+|\gamma(s_\circ)|)^3 [\omega]_3$.  
We also denote $\mathrm M'(x,t)=(\mathrm M_{\gamma}^1 (s_\circ)x,t)$. Since $\mathrm M_{\gamma}^\delta(\scc)=\mathrm M_{\gamma}^1 (s_\circ)
\mathrm D_\delta$,  we have  $\widetilde \omega=\delta^{3} (\overline\omega_{\mathrm M'})^\delta 
$ (see Lemma \ref{weight} for its definition). Using \eqref{scaling-weight} and \eqref{scaling-weight2}, we get  $[\widetilde \omega]_3\le C |\det\mathrm M_{\gamma}^1 (s_\circ)|^{-1} (1+ \|\mathrm M_{\gamma}^1 (s_\circ)\|)^3 [\overline \omega]_3$ because $\det \mathrm M'=\det\mathrm M_{\gamma}^1 (s_\circ)$ and $\|\mathrm M'\|\le 1+ \|\mathrm M_{\gamma}^1 (s_\circ)\|$. 
Combining these two inequalities gives \eqref{weight-}. 
 
 To complete the proof it remains to show \eqref{scale3}.  Changing variables $s \rightarrow \delta s+\scc$  and using  \eqref{scaledcurve}, we see  
$
 A^\gamma[\psi_{\!J}] f(x,t)=\delta \int f \big( x- t\gamma(\scc)- t\mathrm M_{\gamma}^\delta (\scc)\gamma^\delta_{\scc}(s) \big) \psi_{\!J}(\delta s+\scc) ds. 
$    We thus  have 
 \begin{align*}
 A^\gamma[\psi_{\!J}] f(x,t)&=\delta   |\!\det \mathrm M_{\gamma}^\delta(\scc)|^{-\frac1p}\!\!
\int \widetilde f \big( (\mathrm M_{\gamma}^\delta (\scc))^{-1}(x-t\gamma(\scc)) - t\gamma^\delta_{\scc}(s) \big) \psi_{\!J_\circ}(s)\,ds.
\end{align*}
 Therefore the change of variables $x\to \mathrm M_{\gamma}^\delta(\scc) x+t\gamma(\scc)$ yields
\eqref{scale3}. 
\end{proof}

\subsubsection*{Reduction} Let $\gamma\in \mathrm C^D(\mathbb J)$ be a curve satisfying \eqref{nonv}. For a given $\varepsilon_\circ>0$ we take  $\delta= \delta_\ast$ where  $\delta_\ast$ is the number given in Lemma \ref{lem:stable}.   
Applying \eqref{Jsum}  to $\psi \in \mathrm C_0^D(\mathbb J)$ and then  Lemma  \ref{lem:scaled}  to each interval $J$, we   have 
\[  \|  A^\gamma[\psi] f\|_{L^p(\mathbb R^3 \times I, \omega)} \le  \delta^{1-\frac{3}p}  \sum_{J\in \mathfrak J(\delta)}
\big \| A^{ \gamma_{c_J}^\delta} [ \psi^J] \widetilde f^J\big\|_{L^p(\mathbb R^3 \times I, \widetilde \omega^J)},
 \] 
 where $ \gamma_{c_J}^\delta\in \mathfrak C^D(\vepc)$ (by Lemma \ref{lem:stable}),  $[ \widetilde \omega^J]_3\le  C_J[\omega]_3$,  
 $C^{-1}\psi^J\in\mathfrak N^D(J_\circ)$ for some  constants  $C_J,$ $C>0$, and   $\widetilde f^J$  satisfies that $\|\widetilde f^J\|_p\le \|f\|_p$ and 
$\supp \mathcal F(\widetilde  f^J)\subset \{\xi: (B^J)^{-1}\lambda\le |\xi|\le B^J\lambda\}$ for a constant $B^J$. 
Since there are at most  $C\delta_\ast^{-1}$ many intervals, for the estimate  \eqref{Weighted}  it is enough to obtain estimate for each $A^{ \gamma_{c_J}^\delta} [ \psi^J] \widetilde f^J$ against the weight $\widetilde \omega^J$. Hence, in order to show \eqref{Weighted},  after decomposing  $\widetilde f^J$ via Littlewood-Paley projection and replacing $\widetilde \omega^J$ with $(C_J [\omega]_3)^{-1}\widetilde \omega^J$,   we need only to consider 
the curve $\gamma\in \mathfrak C^D(\vepc)$ and the weight $\omega$ with  $[\omega]_3\le 1$.  

Furthermore,  since $A^\gamma[\psi]  f(x,t)=A^\gamma[\psi]  f(\frac{\cdot}r)(rx,rt)$, by scaling
after splitting  $I$ into three intervals $[2^{-1}, 1], [1,2]$ and $[2,4]$,   the proof of   \eqref{Weighted}  now reduces to showing 
\[ \| A^\gamma[\psi]   f\|_{L^p(\mathbb R^3 \times [1,2], \omega)}
\le C \lambda^{-\varepsilon_p} \|f\|_{L^p(\mathbb R^3)}
\] 
 for $[ \omega ]_3\le 1$, $\gamma\in \mathfrak C^D(\vepc)$, and $\psi\in \mathfrak N^D(J_\circ)$ for some $D$.

\begin{defn}
Fixing $p,\varepsilon_\circ, D$,  for $\lambda\ge 1$ we  define the quantity $Q(\lambda)$
by 
\begin{align*}
Q(\lambda)=\sup \big\{ \| 
A^{\gamma} [\psi] f &\|_{L^p(\mathbb R^3\times [1,2],\omega)} :
\gamma \in \mathfrak C^D(\varepsilon_\circ), \,  \psi\in  \mathfrak N^D(J_\circ),  
\\
&\qquad\qquad  [\omega]_3\le 1,\, \supp \widehat f \subset \mathbb A_\lambda, \,
\|f\|_{L^p(\mathbb R^3)}\le 1  \big\}.
\end{align*}  
An elementary estimate  gives $Q(\lambda) \le  C \lambda^2$ for $1\le p\le \infty$.
\end{defn}

Thanks to the discussion in the above and Lemma \ref{WtoM},  Theorem \ref{LSA} now follows from 
the next proposition, which we prove in Section \ref{s-conclusion}.     
 
 \begin{prop} For $p\in (3,\infty)$,   there are positive  constants  $\vepc$,  $D$, $\vep_p$, and $C$  such that 
 \label{QQ}
 \begin{align}\label{Q}
Q(\lambda) \le C   \lambda^{-\varepsilon_p}.
\end{align}
\end{prop}

In order to show $\eqref{Q}$ we need only to handle  $A^\gamma[\psi]$ with $\psi\in  \mathfrak N^D(J_\circ)$,  which  we decompose  in the fashion of  \eqref{Jsum}. 
Thus it suffices  work with the  intervals $J\cap J_\circ\neq \emptyset 
$.  We  set 
\[
\mathfrak J_\circ(\delta)=\big\{J\in\mathfrak J(\delta): J\subset (1+2\cc)J_\circ\big\}.
\]
What follows next  is a   consequence of  Lemma \ref{lem:scaled}, which plays an important role in proving \eqref{Q}. 

\begin{lem}\label{scaling2}
Let $J\in  \mathfrak J_\circ(\delta)$  and  $\psi_{\!J}\in \mathfrak N^D(J)$. Suppose  $\gamma\in \mathfrak C^D(\vepc)$, $[\omega]_3\le 1$,  and   $\supp \widehat  f\subset \mathbb A_\lambda$. If $\delta^3\lambda\ge 2^2$ and $\varepsilon_\circ>0$ is sufficiently small,  there is a constant  $C$, independent of $\gamma,$ $\omega,$ and $\psi_{\!J}$, such that 
\Be
\label{hoho}
\big\|  A^{\gamma}[\psi_{\!J}] f \big\|_{L^p(\mathbb R^3\times [1,2],\omega)}
\le C
\delta^{1-\frac {3}p}  
K_\delta(\lambda) \|f\|_{L^p(\mathbb R^3)},
\Ee
where 
\[K_\delta(\lambda)=\sum_{2^{-2}\delta^3 \lambda \le 2^{j} \le 2^2\delta \lambda} Q(2^{j}).\] 
\end{lem}

\begin{proof}  
We denote  $J=[\scc-\cc\delta, \scc+\cc\delta]$. Since $\gamma\in \mathfrak C^D(\vepc)$, $\gamma_{\scc}^\delta\in \classe$  by Lemma \ref{lem:stable} and our choice of $\delta$, i.e., \eqref{del-choice}. 
Noting that $s_\circ \in 2J_{\circ}$, $\gamma\in \mathfrak C^D(\vepc)$ and $\vepc\le c_\circ^2$,  we see that   $|\gamma(\scc)|\le 3\cc$ and 
$\| \mathrm M_{\gamma}^1(s_\circ)-\mathrm I_3\|\le 5\cc$. 
If we use $ \sum_{\ell=0}^\infty (\rm I_3- \mathrm M_{\gamma}^1(s_\circ))^\ell=(\mathrm M_{\gamma}^1(s_\circ))^{-1}$,  it follows $ \| (\mathrm M_{\gamma}^1 (s_\circ))^{-1}-\mathrm I_3\|\le \frac{5\cc}{1-5\cc}$. Since $\|\mathrm M\|=\|\mathrm M^t\|$ for any matrix $\mathrm M$, $ \| (\mathrm M_{\gamma}^1 (s_\circ))^{-t}-\mathrm I_3\|< 1/100$. 
So, we have
\[ \frac{99}{100}\le  \inf_{|z|=1} |(\mathrm M_{\gamma}^1 (s_\circ))^{-t} z|, \quad
  \|(\mathrm M_{\gamma}^1 (s_\circ))^{-t}\|, \quad  |\det\mathrm M_{\gamma}^1 (s_\circ)|\le \frac{101}{100}.\]
Therefore, by \eqref{weight-} and \eqref{supp-} we  see,  respectively,  that   $[\widetilde \omega]_3\le C$ with a constant $C$ independent of $\gamma$ and  that  
$\supp \mathcal F(\widetilde f\,)\subset  \{ \xi:  2^{-1} \delta^3\lambda \le |\xi |\le 2\delta \lambda\}$.
 Let 
$\beta_\ast\in \mathrm C_0^\infty([3/4,7/4])$  be such that   $\sum_j \beta_\ast (2^{-j}\cdot)=1$. 
We decompose  
\[ \widetilde f= \sum_{ 2^{-2}\delta^3 \lambda \le 2^{j} \le 2^2\delta \lambda } \widetilde f_j,\]
where $\widetilde f_j =\mathcal F^{-1}\big(\beta_\ast(2^{-j}|\cdot|) \mathcal F(\widetilde f\,)\big)$. 
By  \eqref{scale3}  it follows that 
\[  \big\|  A^{\gamma}[\psi_{\!J}] f \big\|_{L^p(\mathbb R^3\times [1,2],\omega)}\le  \delta^{1-\frac{3}p} 
\sum_{2^{-2}\delta^3 \lambda \le 2^{j} \le 2^2\delta \lambda} 
 \| A^{\gamma_{\scc}^\delta} [ \psi_{\!J_\circ}] \widetilde f_j\|_{L^p(\mathbb R^3 \times [1,2], \widetilde \omega)} .
\]
Since $\supp \mathcal F(\widetilde {f_j}) \subset \mathbb A_{2^{j}}$  and $\| \widetilde f_{j} \|_p \le C_{\beta_\ast}\|f\|_p$ and since  $\gamma_{\scc}^\delta\in \mathfrak C^D(\vepc)$, $\psi_{\!J_{\circ}} \in \mathfrak N^D(J_{\circ})$ and $[\widetilde \omega]_3\le C$, we have
$\| A^{\gamma_{\scc}^\delta} [ \psi_{\!J_\circ}] \widetilde f_j\|_{L^p(\mathbb R^3 \times [1,2], \widetilde \omega)}\le CQ(2^{j})\|f\|_p$ while $C$ is  independent of $\gamma,$ $\omega,$ and $\psi_{\!J}$.
Therefore we get \eqref{hoho}.
\end{proof}

\subsection{Decomposition on the  Fourier side}
To show the inequality  \eqref{Q} we need only to deal with $\gamma\in \classe$ and $\psi\in  \mathfrak N^D(J_\circ)$, therefore  it suffices to consider the curve 
$\gamma$ over the interval $(1+2\cc) J_\circ$.  This additional localization helps to simplify the argument, which comes after.  

Since $\vepc\le c_\circ^2$, 
it is clear that 
\Be
\label{e123}
|\gamma'(s)-e_1|\le 2\cc, \  |\gamma''(s)-e_2| \le 2 \cc, \ |\gamma'''(s)-e_3| \le 2 \cc
\Ee 
for $s\in  (1+2\cc)J_\circ$ and $\gamma\in \classe$.  Thus we have $|\gamma'(s)\cdot \xi|  +|\gamma''(s)\cdot \xi| \ge \cc |\xi|$  if $|\xi_1|\ge 3\cc |\xi|$ or  $|\xi_2|\ge 3\cc |\xi|$.  
 Using  Proposition \ref{JC*}  below we can handle the contribution from the  part of frequency $|\xi_1|\ge 3\cc |\xi|$ or  $|\xi_2|\ge 3\cc |\xi|$ since the condition \eqref{nondegenerate} is satisfied.   We shall mainly  concentrate on the case  where $\xi$ is included in the set 
\[
\mathbb A^\ast_\lambda:=\big\{  \xi: 2^{-1}\lambda \le   |\xi|\le 2\lambda, \  |\xi_1|\le  2^{2}\cc |\xi|, \ |\xi_2|\le 2^2\cc |\xi|\big \}.
\]
The following is easy to see.  
 
 \begin{lem} 
 \label{sigma} There exists a function $\sigma \in \mathrm C^{D-2}(\mathbb A^\ast_\lambda) $,  homogeneous of degree $0$,
such that,   for $\xi\in \mathbb A^\ast_\lambda$,  $|\sigma(\xi)|\le 5\cc$ and 
\[
\gamma''(\sigma(\xi))\cdot \xi=0. 
\]
 \end{lem} 
\noindent Indeed, we need to solve the equation $\gamma''(s)\cdot \xi=0$ for a given $\xi$,  equivalently,  ${\xi_3^{-1}}{\xi_2}+ s  +  e(\xi, s)=0$ where $e(\xi, s)$ is a function of 
homogeneous of degree zero and $\|e(\xi, \cdot) \|_{\mathrm C^{D-2}}\le 2\vepc$.  An elementary argument shows existence of $\sigma(\xi)$ and the implicit function theorem guarantees that  $\sigma\in \mathrm C^{D-2}(\mathbb A^\ast_\lambda)$ since $\gamma\in \classe$.   It is clear that   $|\sigma(\xi)|\le 5\cc$
because ${\xi_3^{-1}}{\xi_2}+ \sigma(\xi)  +  e(\xi, \sigma(\xi))=0$.

For  $\xi \in \mathbb A^\ast_\lambda$, we  denote 
\begin{align*}
\Lambda_\gamma(\xi)&= \gamma'''(\sigma(\xi)) \cdot \xi,\\
 R_\gamma(\xi)&= - \frac{\gamma'(\sigma(\xi))\cdot \xi}{\Lambda_\gamma(\xi)}.
\end{align*}
If  $\xi\in \mathbb A^\ast_\lambda$ and $\sigma(\xi)\in  (1+2\cc) J_\circ$,  by \eqref{e123} we have  $2^{-2} \lambda\le |\Lambda_\gamma(\xi)|\le 2^2 \lambda$,  
$ |\gamma'(\sigma(\xi))\cdot \xi-\xi_1| \le 2^3 \cc\lambda $, and 
$\big |\Lambda_\gamma(\xi) -\xi_3\big|\le 2^3  \cc\lambda$, so $|R_\gamma(\xi)|\le 2^6 \cc $. 

\subsubsection*{Decomposition of the operator $ A^\gamma[\psi_{\!J}]$}    By a Taylor  expansion we have   
\begin{align}
\label{SJ}
\gamma'(s)\cdot \xi
&=-\Lambda_\gamma (\xi)R_\gamma(\xi) + 2^{-1}\Lambda_\gamma(\xi)(s-\sigma(\xi))^2 +\mathcal O(\varepsilon_\circ \lambda |s-\sigma(\xi)|^3), 
\\
\label{SJ2}
\gamma''(s)\cdot \xi
&=\Lambda_\gamma(\xi) (s-\sigma(\xi)) +\mathcal O(\varepsilon_\circ \lambda |s-\sigma(\xi)|^2)
\end{align}
for $s\in J$ and $\xi\in \mathbb A^\ast_\lambda$. 
Thus  $\gamma'(s)\cdot \xi$ and $\gamma''(s)\cdot \xi$ have  lower bounds if $\sigma(\xi)$ is distanced from  $J$, so it is not difficult to have control over 
 the  contribution from the associated frequency. However, if $\sigma(\xi)$ is close  to $J$ for $\xi \in \supp \widehat f$, the behavior of  $A^\gamma[\psi_J] f$  becomes less favorable.  
  This leads us to define,  for   $K\ge 1$  and  $J\in\mathfrak J_\circ(\delta)$,     
\[
\mathcal R_{\!J}(K)=\big\{\xi: |\gamma'(c_J)\cdot \xi | \le Kc_\circ^2\delta^2 \lambda ,~ |\gamma''(c_J)\cdot \xi |  \le K\cc\delta \lambda,~ 2^{-2}\lambda\le |\xi_3|\le 2^2\lambda\big\}, 
\]
which  contains the unfavorable frequency part of $A^\gamma[\psi_{\!J}]f$.  
Concerning the sets $\mathcal R_{\!J}(K)$ we have the next lemma, which we use later. 

\begin{lem}
Let $\gamma \in \mathfrak C^D(\varepsilon_\circ)$.
If $\varepsilon_\circ>0$ is sufficiently small,  we have 
the following with  $C$ independent of $\gamma$ and $\delta$: 
\begin{align}\label{BRB1}
\sum_{J\in \mathfrak J_\circ(\delta)} \chi_{
\mathcal R_{\!J}(2^6)} \le C.
\end{align}
\end{lem}

\begin{proof} In order to show \eqref{BRB1}  it is sufficient to verify  that  the sets 
$ {\mathbf r}_{\!J}:=\{\xi: \lambda\xi\in  \mathcal R_{\!J}(2^6) \}$ overlap each other at most $C$ many times.  Note that 
$ {\mathbf r}_{\!J}$  is contained in $2^8\cc\delta$ neighborhood of the line $L_J$ passing through the origin with its direction parallel to $\gamma'(c_J)\times \gamma''(c_J)$.   Since ${\mathbf r}_{\!J}\subset \{\xi: 2^{-4}\le |\xi|\le 2^4\}$, it is sufficient to show that  the directions of the lines $L_J$ are separated from each other by  a distance at least $ 2^{-1} \cc\delta$. This in turn follows if we show 
\[  \frac d{ds}\Big( \gamma'(s)\times \gamma''(s)\Big)=-e_2+ \os (5\cc)  \]
for $\gamma \in \classe$ because the distance between the centers $c_J$  of $J$ is at least $ \cc\delta$.
Since $ (d/{ds})\big( \gamma'(s)\times \gamma''(s)\big)=\gamma'(s)\times \gamma'''(s)$, it is enough to show $\gamma'(s)\times \gamma'''(s)=-e_2+ \os (5\cc)$. 
Since  $s\in [-2\cc,  2\cc] $ and $\gamma\in \classe$, $|\gamma'(s)-e_1|\le 2\cc(1+2\cc)$ and $|\gamma'''(s)-e_3|\le \cc^2$. 
Thus, we have  $\gamma'(s)\times \gamma'''(s)=-e_2+ \os (5\cc)$. 
\end{proof}

Let  $\widetilde \beta \in \mathrm C_0^\infty ([2^{-2}, 2^2])$ be such that $\widetilde\beta=1$ on $[2^{-1}, 2]$. 
Then we  set 
\[  {\widetilde \chi}_{\mathcal R_{\!J}}(\xi)= \beta_0\Big(\frac{|\gamma'(c_J)\cdot \xi |}{ 2^5c_\circ^2\delta^2\lambda}\Big) \beta_0\Big(\frac{|\gamma''(c_J)\cdot \xi | }{ 2^5\cc\delta\lambda}\Big)
\widetilde\beta\Big(\frac{|\xi_3| }{\lambda}\Big),\] 
so that  ${\widetilde \chi}_{\mathcal R_{\!J}}$ is supported in $\mathcal R_{\!J}(2^6)$  and 
${\widetilde \chi}_{\mathcal R_{\!J}}(\xi)=1$ if  $\xi\in \mathcal R_{\!J}(2^5)\cap \mathbb A^\ast_\lambda$. We also set
\[
{P_{\!J} f}=\mathcal F^{-1} ({\widetilde \chi}_{\mathcal R_{\!J}} \widehat f \,).
\]
The following is a consequence of  \eqref{BRB1}. 

\begin{lem} 
\label{ortho} If $\vepc$ is small enough,   we have  
$\big(\sum_{J\in \mathfrak J_\circ(\delta)} \|P_{\!J}f\|_p^p\big)^{1/p} \le C \|f\|_p$
 for $2\le p\le \infty$  whenever  $\gamma\in \classe$.
\end{lem}

\noindent The inequality follows by interpolation between the cases $p=2$ and $p=\infty$.   Plancherel's theorem and \eqref{BRB1} give $(\sum_J \|P_{\!J}f\|_2^2 )^{1/2} \le C \|f\|_2$ and the estimate $\max_{J}  \|P_{\!J}f\|_\infty \le C\|f\|_\infty$ is obvious.

\subsubsection*{Decomposition away from the conic surface $\mathcal C_\lambda$}
We further decompose   $A^\gamma[\psi_{\!J}]P_{\!J} f$ on the  Fourier side taking into account  how close $\xi$ is to the conic set  $\mathcal C_\lambda:=\{\xi\in \mathbb A^\ast_\lambda: R_\gamma(\xi)=0\}$.  
To this end we set  \[\widetilde \chi_{\mathbb A^\ast_\lambda}(\xi)=\beta_0\big(\frac{\xi_1}{2\cc|\xi|}\big) \beta_0\big(\frac{\xi_2}{2\cc|\xi|}\big)   \beta(\lambda^{-1}|\xi|).\]
For  $0<\nu\ll 1$, we  define the cutoff functions  $\pi_{\mathbf c}$,  $\pi_{\mathbf e}$, $\pi_{\mathbf o}^1$, and $\pi_{\mathbf o}^0$ by  
\begin{align*}
\pi_{\mathbf c}(\xi)  
    &=
        \widetilde \chi_{\mathbb A^\ast_\lambda}(\xi) \beta_0(\lambda^{\frac 23-2\nu} |R_\gamma(\xi)|),
             \\
\pi_{\mathbf e}(\xi)  
   &=
          \beta(\lambda^{-1}|\xi|) - \widetilde \chi_{\mathbb A^\ast_\lambda}(\xi) \beta_0(\delta^{-100} |R_\gamma(\xi)|),  
\end{align*}
and, for $j=0,1$, 
\begin{align*}
\pi_{\mathbf o}^j(\xi) 
   &=
     \widetilde \chi_{\mathbb A^\ast_\lambda}(\xi) \chi_{\{\xi:  (-1)^{j+1} R_\gamma(\xi) >0\}}(\xi) \big( \beta_0(  \delta^{-100} |R_\gamma(\xi)|)- \beta_0( \lambda^{\frac 23-2\nu}|R_\gamma(\xi)|) \big).
\end{align*}

 The support of $\widetilde \chi_{\mathbb A^\ast_\lambda}$ is contained in $\mathbb A^\ast_\lambda$ and  $\pi_{\mathbf c}+\pi_{\mathbf o}^1+\pi_{\mathbf o}^0+\pi_{\mathbf e}=\beta(\lambda^{-1}|\cdot|)$ almost everywhere.  The functions $\pi_{\mathbf c}$, $\pi_{\mathbf o}^1+\pi_{\mathbf o}^0$, and $\pi_{\mathbf e}$ 
  roughly  split 
the set $ \mathbb A^\ast_\lambda$ into three regions $\{\xi: |R_\gamma(\xi)| \le C \lambda^{2\nu-2/3}\}$,  $\{\xi: C\lambda^{2\nu-2/3}\le  |R_\gamma(\xi)| \le C_1  \delta^{100} \} $, and $\{\xi:  C_1 \delta^{100}\le   |R_\gamma(\xi)|  \}$.  The   division between the first set and the other two   reflects different  asymptotic behaviors of  
the multiplier $ A^\gamma[\psi_{\!J}] (e^{i (\cdot) \cdot \xi})(0,t)$ as $|\xi|\to \infty$. 
The further division of the second and the third sets  is necessitated  by the  transversality condition for the multilinear estimate, which is to be discussed in the next section.

\medskip

We also define  the associated  multiplier operators  ${\mathcal P_{\!\mathbf {c}}}$, $\mathcal P_{\!\mathbf {o}}^1$, $\mathcal P_{\!\mathbf {o}}^0$, and ${\mathcal P_{\!\mathbf {e}}}$ by 
\begin{align*}
\widehat{{\mathcal P_{\!\mathbf {c}}} g}(\xi)=  \pi_{\mathbf c}(\xi)   \widehat g(\xi), 
\  \ 
\mathcal F(\mathcal P_{\!\mathbf {o}}^j g)(\xi) = \pi_{\mathbf o}^j (\xi)     \widehat g(\xi), \ j=0,1, 
\ \
\widehat{{\mathcal P_{\!\mathbf {e}}} g}(\xi)= \pi_{\mathbf e}(\xi)  \widehat g(\xi).  
\end{align*}
Besides, we  set ${\mathcal P_{\!\mathbf {n}}}={\mathcal P_{\!\mathbf {c}}}+\mathcal P_{\!\mathbf {o}}^1 + \mathcal P_{\!\mathbf {o}}^0$.\footnote{The subscripts $\mathbf c$, $\mathbf e$,
stand for the (main) conic region, the (minor) error parts, respectively, while $\mathbf o$ and 
$\mathbf n$ represent   {\it outside of} and  {\it near}  the conic region, respectively.}
Then easy estimates for the kernels  of the operators give   
\begin{align}
\label{easy1}
 \| {\mathcal P_{\!\mathbf {c}}} \|_{p\to p} \le C_1  \lambda^{C},  \ \  \  & \| \mathcal P_{\!\mathbf {o}}^j \|_{p\to p} \le C_1 \lambda^{C}, \ \ j=0,1, 
\end{align}
for $1\le p\le \infty$ and some constants $C, C_1>0$. It is possible to get better bounds if we use the decoupling or the square function estimate  for the cone (for example, \cite{LV, GWZ}) 
but we do not attempt to do so since it is irrelevant to our purpose.  Similarly,   we also have 
\Be
\label{easy2}
 \| {\mathcal P_{\!\mathbf {e}}}  \|_{p\to p}\le C_1  \delta^{-C},  \ \  \| {\mathcal P_{\!\mathbf {n}}} \|_{p\to p} \le C_1  \delta^{-C}
\Ee 
for $1\le p\le \infty$.  
For the former we need only to note that $\|  \mathcal F^{-1}(\pi_{\mathbf e}) \|_{L^1(\mathbb R^3)}\le C_1 \delta^{-C}$. The latter follows from the former   because  the  multiplier associated to the operator   ${\mathcal P_{\mathbf n}}$ is $\beta(\lambda^{-1}|\cdot|)- \pi_{\mathbf e}$.

\subsection{Nondegenerate part}
Decomposition of the operator $A$ on the Fourier side gives rise to operators of  the form of \eqref{Aaf} such as  $A^\gamma[\psi_{\!J}] P_{\!J}$, $A^\gamma[\psi_{\!J}] (1-P_{\!J})$, \dots, $A^\gamma[\psi_{\!J}] \mathcal P_{\mathbf e} $. 
If  $|\gamma'(s)\cdot \xi| + |\gamma''(s)\cdot \xi| \ge C |\xi|$ on the support of $a$, we can handle $A^\gamma[a]$ using  the following  theorem which  is  a straightforward consequence of  \cite[Theorem 4.1]{PS}. 
 
\begin{thm}   
\label{PS6}
Let $K\ge 1$ and  $[s_\circ-2r, s_\circ+2r]\subset {\mathbb J}$ with $K^{-1}\le r$. Suppose that $a(s,t,\xi)$ is a smooth function supported in  $[s_\circ-r, s_\circ+r] \times I \times \mathbb A_\lambda$ and  
$
|\partial_s^{j_1} \partial_t^{j_2} \partial_\xi^\alpha a(s,t,\xi)| \le B |\xi|^{-|\alpha|}
$
 for $|\alpha|\le 5$ and $ j_1, j_2=0,1$. 
Also, assume that 
\begin{align}\label{nondegenerate}
|\gamma'(s)\cdot \xi| + |\gamma''(s)\cdot \xi| \ge K^{-1} |\xi|
\end{align}
holds whenever  $(s,t,\xi) \in \supp a$ for some $t\in I$. 
Then, if $p \ge 6$ and $\vepc>0$ is small enough,  for  any $\vep>0$, 
\begin{align}\label{ps41*}
\|  A^\gamma[a] f \|_{L^p(\mathbb R^3\times I)} 
\le 
C_\varepsilon B K^{C}  \lambda^{-\frac 2p+\varepsilon} \|f\|_{L^p(\mathbb R^3)}
\end{align}
holds  whenever $\gamma\in \mathfrak C^D(\vepc)$ and $\widehat f$ is  supported in $\mathbb A_\lambda$. 
\end{thm}

The statement of Theorem \ref{PS6}  differs from the one  in \cite{PS} in a couple of aspects. First, the range of $p$ is enlarged to $p\ge 6$ \footnote{The critical case $p=6$ can be included by interpolation with a trivial estimate.} thanks to  the $\ell^p$-decoupling inequality for the cone \cite{BD}. 
Secondly, there is  an extra factor $K^{C}$  in \eqref{ps41*}.  One can easily show the estimate \eqref{ps41*}  by following the argument in  \cite{PS}.
It is also possible to deduce  \eqref{ps41*} from that with $K\in [2^{-1},2]$ by finite decomposition  and making use of scaling and affine transform.  Uniformity of the bound over $\gamma\in \mathfrak C^D(\vepc)$ is  clear.

The estimate $|\int e^{-i t\gamma(s)\cdot \xi}  a(s,t,\xi) ds|\le C_1BK^C|\xi|^{-\frac12}$ follows by \eqref{nondegenerate} and van der Corput's Lemma. We thus have 
$\| A^\gamma[a] f \|_{L^2(\mathbb R^3\times I)} 
\le C_1B K^{C}  \lambda^{-\frac12} \|f\|_{L^2(\mathbb R^3)}$  by Plancherel's theorem.  Interpolation between the estimate  and \eqref{ps41*} with $p=6$ gives

\begin{cor}
\label{cor:nondegenerate}  
Under the same assumption as in Theorem \ref{PS6},  suppose $2\le p \le 6$ and $\vepc$ is small enough. Then, for  any $\vep>0$, 
\begin{align*}
\|  A^\gamma[a] f \|_{L^p(\mathbb R^3\times I)} 
\le 
C_\varepsilon B K^{C} \lambda^{-\frac 14-\frac1{2p}+\varepsilon} \|f\|_{L^p(\mathbb R^3)}
\end{align*}
holds whenever $\gamma\in \mathfrak C^D(\vepc)$ and $\widehat f$ is  supported in $\mathbb A_\lambda$. 
\end{cor}

We also make use of the following (\cite[Theorem 1.4]{PS}).

\begin{thm}[]
Let $J\subset {\mathbb J}$ be a compact interval of length $\delta$ and $\psi_{\!J}\in {\mathfrak N^D(J})$. 
Then, if $p\ge6$ and $\vepc$ is small enough,  
\begin{align}
\label{ps14*}
\|  A^\gamma[\psi_{\!J}] f \|_{L^p(\mathbb R^3 \times I )} \le C_\varepsilon \delta^{-C} 
\lambda^{-\frac 4{3p}+\varepsilon} \|f\|_{L^p(\mathbb R^3)}
\end{align}
holds  for any $\vep>0$ whenever  $\gamma \in \classe$ and  $\widehat f$ is supported in $\mathbb A_\lambda$. 
\end{thm}

\noindent Compared with  \cite[Theorem 1.4]{PS},   the range of $p$ is extended to $p\ge 6$ by the aforementioned decoupling inequality \cite{BD}. 
The estimate \eqref{ps14*} with additional factor $\delta^{-C}$ can be shown by scaling  and  its uniformity  over $\gamma\in \mathfrak C^D(\vepc)$ is also obvious.

\subsubsection*{Estimates for  $A^\gamma[\psi_{\!J}] (1-P_{\!J})$  and  $ A^\gamma[\psi_{\!J}] \mathcal P_{\mathbf e}$} 
The  condition \eqref{nondegenerate} is satisfied on the support of  $\psi_{\!J}(s)(1-{\widetilde \chi}_{\mathcal R_{\!J}}(\xi))$. Thus, 
 using  Corollary   \ref{cor:nondegenerate}, we can get a favorable estimate for $A^\gamma[\psi_{\!J}] (1-P_{\!J})$. 
 We also obtain a similar estimate for $ A^\gamma[\psi_{\!J}] \mathcal P_{\mathbf e}$ (see Proposition \ref{JC*} below).

\begin{prop}\label{JC}
Let  $[\omega]_3\le 1$, and 
 $J\in \mathfrak J_\circ(\delta)$. 
 If  $2\le p \le 6$ and $\vepc>0$ is small enough,  
\begin{align}\label{small}
\|  A^\gamma[\psi_{\!J}] (1-P_{\!J})f\|_{L^p(\mathbb R^3 \times [1,2], \omega)}
\le C_\varepsilon \delta^{-C}  \lambda^{\frac 12(\frac1p-\frac12)+\varepsilon}\|f\|_{L^p(\mathbb R^3)}
\end{align}
holds for  any $\vep>0$  whenever  $\supp \widehat f \subset \mathbb A_\lambda$, $\gamma \in \mathfrak C^D(\varepsilon_\circ)$,  and $\psi_{\!J}\in \mathfrak N^D(J)$. 
\end{prop}

We prove Proposition \ref{JC} using the next lemma. To this end, we need only to have the inequality \eqref{control}  below for $1\le p\le \infty$.  However, for later use (see Section  \ref{sec:3.4}) we prove it for $0<p\le \infty$.

\begin{lem}\label{change}
Let $0<p\le \infty$, $0<\alpha \le 4$ and $\omega\in \Omega^\alpha$. Suppose that $F \in L^p(\mathbb R^4)$ and $\widehat F$ is supported on $\mathbb B^4(0,\lambda)$. Then we have
\begin{align}\label{control}
\|F\|_{L^p(\mathbb R^4,\omega)}
\le C
[\omega]_\alpha^{\frac 1p} \lambda^{\frac {4-\alpha}p}\|F\|_{L^p(\mathbb R^4)}.
\end{align}
\end{lem}

\begin{proof}  Since $wdxdt\ll dxdt$ for $\omega\in \Omega^\alpha$,   it follows that $\|F\|_{L^\infty(\mathbb R^4,\omega)}\le \|F\|_{L^\infty(\mathbb R^4)}$. 
When $1\le p<\infty$,  \eqref{control} is a simple consequence of H\"older's inequality. Indeed, 
let us take $\varphi \in \mathcal S(\mathbb R^4)$ such that $\widehat \varphi=1$ on $\mathbb B^4(0,1)$
and $\widehat \varphi$ is supported on $\mathbb B^4(0,2)$. 
Then $F=F\ast \varphi_\lambda$ since $\widehat F$ is supported on $\mathbb B^4(0,\lambda)$. 
Thus, by H\"older's inequality we have 
$|F|^{p} \le C |F|^{p} \ast  |\varphi|_\lambda$. So, 
\[
\| F\|_{L^p(\mathbb R^4, \omega)}^p
\le C \! \int |F(x)|^p  
|\varphi|_\lambda \ast \omega(x) dx
\le  C\|F\|_{L^p(\mathbb R^4)}^p \| |\varphi|_\lambda \ast\omega \|_\infty.\]
 This gives \eqref{control} since $\||\varphi|_\lambda \ast \omega\|_\infty \le C  \lambda^{4-\alpha} [\omega]_\alpha$.

When $p\in (0,1)$, we  claim that 
\Be\label{convolution}
|F|^{p} \le C |F|^{p} \ast  |\varphi|_\lambda^p,
\Ee
where we denote $|\varphi|_\lambda^p=(|\varphi|^p)_\lambda$. 
Once we have \eqref{convolution}, the proof of \eqref{control} is straightforward. By \eqref{convolution} we have 
$
\| F\|_{L^p(\mathbb R^4, \omega)}^p
\le  C\|F\|_{L^p(\mathbb R^4)}^p \| |\varphi|_\lambda^p \ast\omega \|_\infty.$ 
Since $\||\varphi|_\lambda^p \ast \omega\|_\infty \le C  \lambda^{4-\alpha} [\omega]_\alpha$, we obtain   \eqref{control}.

We now show \eqref{convolution}.   By scaling we may assume $\lambda=1$,  otherwise one may replace $F$ with  $ F(\cdot/\lambda)$. 
To show \eqref{convolution} when $\lambda=1$,  we first notice   that
\[ 
|F| \ast |\varphi|(x)=\int |F(y)\varphi(x-y)| dy
\le 
\| F\varphi(x-\cdot) \|_\infty^{1-p} \,   |F|^p \ast |\varphi|^p(x). 
\] 
Since  the Fourier transform of $F\varphi(x-\cdot)$ is supported in $\mathbb B^4(0,5)$,  $F\varphi(x-\cdot)=\big(F \varphi(x-\cdot) \big)\ast 5^{4}\varphi(5\,\cdot)$ 
and thus 
$  \|F\varphi(x-\cdot) \|_\infty\le C  \|F\varphi(x-\cdot) \|_1= C  |F| \ast |\varphi|(x)$.
Combining this and the inequality above  gives  $(|F| \ast |\varphi|)^p\le C  |F|^p \ast |\varphi|^p$. 
Since $|F| \le  |F| \ast |\varphi|$, we get \eqref{convolution} with $\lambda=1$.
\end{proof}

\begin{proof}[Proof of Proposition \ref{JC}]
We set \[a(s,t,\xi) = \widetilde\chi(t) \psi_{\!J}(s) (1-{\widetilde \chi}_{\mathcal R_{\!J}}(\xi))\beta(\lambda^{-1}|\xi|),\] so that  $ A^\gamma[a]f= \widetilde\chi(t)A^\gamma[\psi_{\!J}] (1-P_{\!J})f$.  We claim that \eqref{nondegenerate} holds on the support of $a$ with  $K=C_1\delta^{-2}>0$,  $C_1=C_1(\cc)$. 

To see this, it suffices to consider the case   $\xi \in \mathbb A^\ast_\lambda$  because of \eqref{e123}.
We first note that   $
\xi \in  \mathcal R_{\!J}(2^5)$ if   $\sigma(\xi) \in  [\,c_J-|J|,  c_J+|J|\,]$ and $\ |R_\gamma(\xi)| \le  2^3 c_{\circ}^2\delta^2$. 
 Indeed, since $|\sigma(\xi)-c_J| \le 2\cc \delta$,  we have $|\gamma'(c_J)\cdot \xi| \le  2^5 c_\circ^2\delta^2\lambda$ by \eqref{SJ}  and   $|\gamma''(c_J)\cdot \xi | \le 2^3\cc\delta\lambda$ by \eqref{SJ2}. So, it follows $
\xi \in  \mathcal R_{\!J}(2^5)$ since 
$\xi\in \mathbb A^\ast_\lambda$.   Thus, if  $\xi\in \supp (1-{\widetilde \chi}_{\mathcal R_{\!J}})\beta(\lambda^{-1}|\cdot|)\cap  \mathbb A^\ast_\lambda$,  
we have $\sigma(\xi) \notin [c_J-|J|,  c_J+|J|\,]$
or $|R_\gamma(\xi)| \ge  2^3c_\circ^2\delta^2$. 
In the first case,   by  \eqref{SJ2}  we see  $|\gamma''(s)\cdot \xi| \ge 2^{-2}  \cc \delta\lambda$ for all $s\in \supp \psi_{\!J}$. So, we may assume $|s-\sigma(\xi)| \le  3 \cc \delta$ and  $|R_\gamma(\xi)| \ge 2^3c_\circ^2\delta^{2}$. Then we get $|\gamma'(s)\cdot \xi| \ge 2  c_\circ^2\delta^{2}\lambda$ using \eqref{SJ}.   This shows the claim. 

Since \eqref{nondegenerate} holds on the support of $a$,  by Corollary \ref{cor:nondegenerate} we have the estimate  
 \begin{align}\label{small1}
\|   \widetilde\chi(t)A^\gamma[\psi_{\!J}] (1-P_{\!J})f\|_{L^p(\mathbb R^3 \times \mathbb R)}
\le C_\varepsilon \delta^{-C}  \lambda^{-\frac14-\frac1{2p}+\varepsilon}\|f\|_{L^p(\mathbb R^3)}
\end{align}
for $2\le p\le 6$.   We use the estimate  to obtain the weighted estimate \eqref{small}. Since the argument is similar to  that in the proof of Lemma  \ref{WtoM}, we shall be brief.

As before, let us define an  operator $\widetilde  A_J$ by 
\[
\mathcal F(\widetilde  A_J  h)(\xi,\tau)
=
\beta_0 ((\lambda r_0)^{-1} \tau)  \beta(\lambda^{-1} |\xi|)
\mathcal F\big(\widetilde \chi(t) A^\gamma [{\psi_{\!J}}] h\big)(\xi,\tau),
\]
where $r_0=  1+ 4\max\{|\gamma(s)|: s \in {\rm\supp} \psi_{\!J}\} $. Then we have  
$|(\widetilde \chi(t) A^\gamma [\psi_{\!J}]  - \widetilde  A_J)h|\le  C\widetilde E_t^N\ast |h|$ for any $N$ if we use Lemma \ref{lem:ker2}. 
Putting together this (e.g., \eqref{At-A}), $[\omega]_3\le 1$ and  $ \|(1-P_{\!J}) f\|_p\le C \|f\|_p$,  
 we see that 
\begin{align*}
\|\widetilde \chi(t) A^\gamma[\psi_{\!J}] (1-P_{\!J}) f\|_{L^p(\mathbb R^3 \times \mathbb R, \omega)}
&\le \| \widetilde  A_J(1-P_{\!J}) f\|_{L^p(\mathbb R^3 \times \mathbb R, \omega)} + C \lambda^{-N}\|f\|_p.
\end{align*}
The Fourier transform of $\widetilde  A_J(1-P_{\!J}) f$ is supported in $\mathbb B^4(0, 2^2r_0\lambda)$. 
By Lemma \ref{change} we have   $\| \widetilde  A_J(1-P_{\!J}) f\|_{L^p(\mathbb R^3 \times \mathbb R, \omega)}\le C  \lambda^{1/p} \| \widetilde  A_J(1-P_{\!J}) f\|_{L^p(\mathbb R^3 \times \mathbb R)}$. 
Disregarding the minor contribution of  $(\widetilde \chi(t) A^\gamma [\psi_{\!J}]  - \widetilde  A_J)(1-P_{\!J}) f$,  
we only need to consider  $\widetilde \chi(t) A^\gamma [\psi_{\!J}] (1-P_{\!J}) f$ in  $L^p(\mathbb R^3 \times \mathbb R)$. 
Therefore we obtain the estimate \eqref{small} by \eqref{small1}. 
\end{proof}

\begin{prop}\label{JC*} 
Under the same assumption as in  Proposition \ref{JC}, if  $2\le p \le 6$ and $\vepc>0$ is small enough,  for any $\vep>0$,  
\begin{align*}
\|  A^\gamma[\psi_{\!J}] \mathcal P_{\mathbf e} f\|_{L^p(\mathbb R^3 \times [1,2], \omega)}
\le C_\varepsilon \delta^{-C}  \lambda^{\frac 12(\frac1p-\frac12)+\varepsilon}\|f\|_{L^p(\mathbb R^3)}
\end{align*}
holds whenever  $\supp \widehat f \subset \mathbb A_\lambda$, $\gamma \in \mathfrak C^D(\varepsilon_\circ)$ and $\psi_{\!J}\in \mathfrak N^D(J)$. 
\end{prop}

\begin{proof}   We set $\pi_{\bf e}^1(\xi)= \widetilde \chi_{\mathbb A^\ast_\lambda}(\xi) \big( 1- \beta_0(\delta^{-100} |R_\gamma(\xi)|) \big) $ and  
$ \pi_{\bf e}^2(\xi)=\beta(\lambda^{-1}|\xi|)-\widetilde \chi_{\mathbb A^\ast_\lambda}(\xi)$, so that $\pi_{\bf e}=\pi_{\bf e}^1+\pi_{\bf e}^2$. Then we break 
$\widetilde \chi(t) A^\gamma[\psi_{\!J}] \mathcal P_{\mathbf e} f= A^\gamma [a^1] f+  A^\gamma [a^2] f$, where 
\[ a^j(s,t, \xi)= \widetilde \chi(t)\psi_{\!J}(s)\pi_{\bf e}^j (\xi),  \ j=1,2.\]

We first  consider  $A^\gamma[a^1]  f$.  After decomposing   $\psi_{\!J}$ into the bump functions $\psi_\ell$ supported in finitely overlapping intervals $J_\ell$ such that  $\delta^{100}\le |J_\ell|\le  2\delta^{100}$,  $\psi_{\!J}=\sum\psi_{\ell}$, and $|\psi^{(k)}_\ell|\le C_k\delta^{-100k}$,  we set $a^1_\ell(s,t,\xi) =\widetilde \chi(t)\psi_{\ell}(s)\pi_{\bf e}^1(\xi)$.
By \eqref{SJ2}  $|\gamma''(s)\cdot\xi|\ge 2^{-3} \lambda\delta^{100}$  for $s\in \supp  \psi_\ell  $ if $\sigma(\xi)\not\in 
[c_{J_\ell}-|J_\ell|, c_{J_\ell}+|J_\ell|]$. Otherwise,  from  \eqref{SJ}  we have $|\gamma'(s)\cdot\xi|\ge 2^{-2} \delta^{100}\lambda$ for $s\in \supp  \psi_\ell $  since  $|R_\gamma(\xi)|\ge \delta^{100}$  on $\supp \pi_{\mathbf e}^1 $. Therefore  \eqref{nondegenerate} holds with $K=C\delta^{-100}$ for $(s,t,\xi)\in  \supp a^1_\ell$. Applying Corollary \ref{cor:nondegenerate}  with $a=a^1_{\ell}$ gives  
\[ \|   A^\gamma[a^1_{\ell}] f  \|_{L^p(\mathbb R^3\times \mathbb R)}  \le 
C_\vep \delta^{-C} \lambda^{-\frac 14-\frac1{2p}+\varepsilon} \|f\|_{L^p(\mathbb R^3)}. \] 
Arguing similarly as in the proof of  Proposition \ref{JC}, we get  the weighted estimate 
$\|  A^\gamma[a^1_{\ell}]  f\|_{L^p(\mathbb R^3 \times [1,2], \omega)}  \le 
C_\vep \delta^{-C}  \lambda^{-\frac 14+\frac1{2p}+\varepsilon} \|f\|_{L^p(\mathbb R^3)}$. 
Summation over $\ell$  thus gives 
the desired estimate  since there are at most $C\delta^{-100}$ many $\ell$.

The estimate $\|  A^\gamma[a^2 ]  f\|_{L^p(\mathbb R^3 \times [1,2], \omega)}$ $ \le 
C_\vep \lambda^{-\frac 14+\frac1{2p}+\varepsilon} \|f\|_{L^p(\mathbb R^3)}$ can be obtained  likewise but more straightforwardly since $|\gamma'(s)\cdot \xi|  +|\gamma''(s)\cdot \xi| \ge  \cc  |\xi|$ on $\supp a^2$.    
\end{proof}

\section{Multilinear estimates}\label{sec3}

The main object of this section is to prove the following weighted multilinear (quadrilinear) estimate for $ A^\gamma[\psi_{\!J} ]  {\mathcal P_{\!\mathbf {n}}}   f$.
Throughout this section we assume $\gamma\in \classe$ with an $\vepc$ small enough.

\begin{prop}\label{global}  
Let  $J_k \in \mathfrak J_\circ(\delta)$, $1\le k\le 4$, and\, $[\omega]_3\le 1$.  Suppose that $\widehat f_1, \dots, \widehat f_4$   are supported in $\mathbb A_\lambda$ and $\dist(J_\ell, J_k)\ge \delta$, $\ell\neq k$. 
If $14/5<p \le 6$,   there are  constants $\varepsilon_p>0$, $D$, and $C_\delta>0$ such that
\begin{align}\label{gb}
\Big\| \prod_{k=1}^4 
| A^\gamma[\psi_{\!J_k}] ( {\mathcal P_{\!\mathbf {n}}}  P_{\!J_k} f_k)|^{\frac 14} 
\Big\|_{L^p(\mathbb R^3 \times [1,2], \omega)}
&\le C_\delta 
\lambda^{-\varepsilon_p} 
\prod_{k=1}^4 \| f_{k}\|_{L^p(\mathbb R^3)}^{\frac 14}
\end{align}
holds whenever $\gamma \in \classe$ and $\psi_{\!J_k}\in \mathfrak N^D(J_k)$,  $1\le k\le 4.$
\end{prop}

\subsection{Expansions of the multiplier}
In order to prove Proposition \ref{global} we first try to express $ A^\gamma[\psi_{\!J_k}] {\mathcal P_{\!\mathbf {n}}}$  as
a sum  of adjoint restriction operators. To do so, we expand the Fourier multiplier  of the operator  $ A^\gamma[\psi_{\!J_k}] {\mathcal P_{\!\mathbf {n}}}$ into a  series of suitable form.
 We  handle  separately $ A^\gamma[\psi_{\!J_k}] {\mathcal P_{\!\mathbf {c}}}$ (Lemma \ref{op2}) and $ A^\gamma[\psi_{\!J_k}] \mathcal P_{\!\mathbf {o}}^j$, $j=1,0$ (Lemma \ref{op3}).
The estimates in Lemma \ref{op2} and \ref{op3} are somewhat  rough, but  they are good enough for our purpose. So we do not attempt to make them as efficient as possible.

\subsubsection*{Multiplier of  $A^\gamma[\psi_{\!J}] {\mathcal P_{\!\mathbf {c}}}$}  Let $J\in \mathfrak J_\circ(\delta)$. For $\psi_{\!J}\in \mathfrak N^D(J)$ we  set 
\[m_{\!J}(t,\xi)= (2\pi)^{-3}\int e^{-it\gamma(s)\cdot \xi} \psi_{\!J} (s)\,ds.\] 
The worst decay in $\xi$   of  the multiplier  $m_{\!J}(t,\cdot)$  is related to the behavior of $\gamma(s)\cdot \xi$ near $s=\sigma(\xi)$. 
We define 
\[
\Phi^{\bf c}(\xi)=\gamma(\sigma(\xi))\cdot \xi, \quad \xi\in  \mathbb A^\ast_\lambda,
\]
and an  adjoint restriction  operator $\mathcal T^{\bf c}_\lambda$  by setting
\[
\mathcal T^{\bf c}_\lambda g(x,t)
=
\int_{\mathcal C_\lambda^{\bf c} (\delta) } e^{i( x\cdot \xi - t\Phi^{\bf c} (\xi) )} g(\xi)\,d\xi,
\]
 where  $\mathcal C_\lambda^{\bf c} (\delta)= \{\xi\in  \mathbb A^\ast_\lambda: | R_\gamma(\xi) |\le 2\delta^{100} \}.$ We note that $\supp \pi_{\bf c}\subset \mathcal C_\lambda^{\bf c} (\delta) $. 

\begin{lem}\label{op2}
Let $0<\nu \ll1$ and  $J\in \mathfrak J_\circ(\delta)$.   
  Suppose  $\gamma \in \mathfrak C^D(\varepsilon_\circ)$, $\psi_{\!J}\in \mathfrak N^D(J)$, and $\widehat f$ is supported on $\mathbb A_\lambda$. 
 Then we have 
\[
 A^\gamma[\psi_{\!J}] {\mathcal P_{\!\mathbf {c}}} f
             =\sum_{{\ell \in \mathbb Z:\, |\ell | \le \lambda^{10\nu}}} e^{it\ell} \mathcal T^{\bf c}_\lambda \big (c_\ell \pi_{\bf c} \widehat  f \,\big)
                             +\mathcal E_{\mathbf c}  f, \quad  t\in I,
\]
and the following hold with  $C$, $C_N$, and  $C_\delta$  independent of $\gamma$ and $\psi_{\!J}$: 
\begin{align}\label{dec-m0}
|c_\ell(\xi)| &\le C_N \lambda^{\nu-\frac 13} (1+\lambda^{-3\nu}|\ell |)^{-N} 
\end{align}
for any $N$ and  
\begin{align}
\label{error-est0}
\|\mathcal E_{\mathbf c}  f\|_{L^q(\mathbb R^3 \times I)} &\le C_\delta \lambda^{C-\frac 32 \nu D}\| f\|_p, \quad  1 \le p\le q \le \infty. 
\end{align}    
\end{lem}

Summation over $\ell$ results from the Fourier series expansion in $t$ of an amplitude function which appears after factoring out $e^{-it\Phi^{\bf c} (\xi)}$. It simplifies the amplitude function depending both on $\xi$ and $t$ which causes considerable loss in its bound when we attempt to directly apply the multilinear restriction estimate (for example, see \cite[Theorem 6.2]{BCT}).

For the proof of Lemma \ref{op2} and \ref{op3} below,  we write ${m_{\!J}}(t,\xi)$ in a different form. 
Changing of variables $s \rightarrow s+\sigma(\xi)$, we have
\Be
\label{mmj}
{m_{\!J}}(t,\xi)
=(2\pi)^{-3} e^{-it\Phi^{\bf c}(\xi)} \int e^{-it\phi(s,\xi)}  \psi_{\!J}(s+\sigma(\xi)) ds,
\Ee
where  
\begin{align*}
\phi(s,\xi)&:=\gamma(s+\sigma(\xi))\cdot \xi- \gamma(\sigma(\xi))\cdot \xi\,. 
\end{align*}
We here note that  $J\subset (1+2\cc)J_\circ$ and   $|\sigma(\xi)|\le 5\cc$ for $\xi\in \mathbb A_\lambda^\ast$ by Lemma \ref{sigma}. 
Thus  $\phi \in \mathrm C^{D-2}([-1/2, 1/2]\times \mathbb A_\lambda^\ast)$ and $\supp \psi_{\!J}(\cdot+\sigma(\xi))\subset 2^3 J_\circ$. 
 Since $\gamma\in \mathfrak C^D(\vep_\circ)$ and $\gamma''(\sigma(\xi))\cdot\xi =0$, by Taylor expansion of $\phi(\cdot, \xi)$ around $s=0$  it follows that 
\begin{align} 
\label{taylor-phi}
&\phi(s,\xi)
=
\Lambda_\gamma(\xi) \big( -R_\gamma(\xi) s+\frac{1}{6} s^3
+\Theta(s,\xi) \big),  
\\     
 \label{taylor-phi2}
& | \partial_s^k \Theta(s,\xi)|\le C_k \varepsilon_\circ  |s|^{\max(4-k, 0)}, \quad 0\le k\le D.
\end{align} 
In what follows we occasionally resort to \eqref{taylor-phi} and \eqref{taylor-phi2} to exploit the properties of  the phase function  $\phi(\cdot,\xi)$.

\begin{proof}[Proof of Lemma \ref{op2}] We need to consider    ${m_{\!J}}(t,\xi)$ while $\xi\in \supp \pi_{\bf c}$.   We break 
\[  \psi_{\!J}(s+\sigma(\xi))= a_m(s,\xi) +  a_e(s,\xi),\] 
where $\displaystyle{
 a_m(s,\xi)= \psi_{\!J}(s+\sigma(\xi)) \beta_0(2^{-4}\lambda^{\frac13-\nu}s). }$
Then we put 
\begin{align*}
\mathcal I_{\theta}(t,\xi) 
&=(2\pi)^{-3} \int e^{-it\phi(s,\xi)}  a_\theta (s,\xi)ds,   \quad \theta\in \{m, e\}.
\end{align*}
By \eqref{mmj}  it follows that 
\[
{m_{\!J}}(t,\xi)
=e^{-it\Phi^{\bf c}(\xi)} \big( \mathcal I_{m}(t,\xi) + \mathcal I_e(t,\xi)\big).
\]
The major term is $\mathcal I_m$ while $\mathcal I_e$ decays fast as $\lambda\to \infty$.    
 Let $\chi_{\circ} \in \mathrm C_0^\infty([0,2\pi])$ 
such that $ \chi_{\circ}=1$ on the interval $[2^{-1}, 2^2]$.  Expanding   $\chi_{\circ} \mathcal I_m(\cdot,\xi)$ into  Fourier series  over the interval $[0, 2\pi]$
we have
\[
\chi_{\circ}(t) {\mathcal I_m}(t,\xi) =\sum_{\ell \in \mathbb Z} c_\ell(\xi) e^{it\ell}. 
\]
Note that 
$\mathcal F(\chi_{\circ} {\mathcal I_m}(\cdot,\xi))(\ell)= (2\pi)^{-3} \int \widehat{\chi_\circ}( \ell+ \phi(s,\xi)) a_m (s,\xi)ds $. Since $|\phi(s, \xi)| \le C \lambda^{3\nu}$ on $\supp a_m(\cdot, \xi)$  by \eqref{taylor-phi},  we have 
$|\mathcal F(\chi_{\circ} {\mathcal I_m}(\cdot,\xi))(\ell)|\le C\lambda^{\nu-\frac 13} |\ell|^{-N}$ for any $N$ if $|\ell|\ge C_1 \lambda^{3\nu} $ for a large $C_1$. Thus
we get  \eqref{dec-m0} for any $N>0$. 
We also note that   $|\partial_\xi^\alpha \phi|\le C$ and $|\partial_\xi^\alpha a_m|\le C_\delta$ because  
$|\partial_\xi^\alpha \sigma|
\le C \lambda^{-|\alpha|} $ on $\mathbb A^\ast_\lambda$ for  $|\alpha|\le D-2$ (see Lemma \ref{sigma}). 
By the same argument we  obtain, for any $N>0$, 
\Be\label{cl} 
|\partial_{\xi}^\alpha c_\ell(\xi)| \le C_\delta \lambda^{\nu-\frac 13} (1+\lambda^{-3\nu}|\ell |)^{-N}.
 \Ee

We  now put 
\[ \mathcal E_{\mathbf c}  g(x,t) =  (2\pi)^{3}  \sum_{|\ell|>\lambda^{10\nu}}  e^{it\ell}  \mathcal F^{-1}_x\big(  c_\ell         e^{-it\Phi^{\bf c}} \pi_{\bf c}\widehat{g} \,\big)  + (2\pi)^{3}   \mathcal F^{-1}_x
\big(\mathcal I_e(t,\cdot) e^{-it\Phi^{\bf c}}  \pi_{\bf c}\widehat{g}\,\big) .\]
We shall  show  \eqref{error-est0} to complete the proof.  The  terms  $\mathcal F^{-1}_x\big(  c_\ell   e^{-it\Phi^{\bf c}} \pi_{\bf c}\widehat{g}\, \big)$ in the sum  can be handled easily.   Combining the estimate  \eqref{cl}  and  $ |\partial_\xi^\alpha e^{-it\Phi^{\bf c}}|\le C$ for $|\alpha|\le 4$,  we see that  $\mathcal F^{-1}_x
(c_\ell         e^{-it\Phi^{\bf c}} \pi_{\bf c}\widehat{g}\, )=K_t\ast  {\mathcal P_{\!\mathbf {c}}} g$ and $|K_t|\le C_\delta  \lambda^C (1+\lambda^{-3\nu}|\ell |)^{-N}(1+|x|)^{-4}$. 
Thus,  the convolution inequality gives 
\[\|\mathcal F^{-1}_x
(c_\ell         e^{-it\Phi^{\bf c}} \pi_{\bf c}\widehat{g}\, )\|_{L^q(\mathbb R^3 \times I)} \le C_\delta   \lambda^C (1+\lambda^{-3\nu}|\ell |)^{-N} \|{\mathcal P_{\!\mathbf {c}}} g\|_p\]
 for $1\le p\le q\le \infty$.  Taking a large $N\ge D$ and using  \eqref{easy1},  we obtain   $\sum_{|\ell|\ge \lambda^{10\nu}}\|\mathcal F^{-1}_x
(c_\ell         e^{-it\Phi^{\bf c}} \pi_{\bf c}\widehat{g} \,)\|_{L^q(\mathbb R^3 \times I)}\le  C_\delta \lambda^{C-2\nu D}\|g\|_p$.

In order to show the estimate for $\mathcal F^{-1}_x 
(\mathcal I_e(t,\cdot)  e^{-it\Phi^{\bf c}}  \pi_{\bf c}\widehat{g}\,)$  we claim 
\Be
\label{Ie}
|\partial_\xi^\alpha \mathcal I_e(t,\xi)|\le C_\delta \lambda^{-\frac 32\nu (D-|\alpha|)}, \quad |\alpha|\le 4
\Ee
for $\xi\in \supp \pi_{\bf c}$. Using \eqref{Ie},  similarly as before,   we  see     $\mathcal F^{-1}_x
(\mathcal I_e(t,\cdot)  e^{-it\Phi^{\bf c}}  \pi_{\bf c}\widehat{g}\,)=K_t\ast  {\mathcal P_{\!\mathbf {c}}} g$ with  $|K_t|\le C_\delta \lambda^{C-\frac 32\nu D}(1+|x|)^{-4}$. 
Therefore,  the convolution inequality and  \eqref{easy1}  give
\[ \big\|\mathcal F^{-1}_x
(\mathcal I_e(t,\cdot) e^{-it\Phi^{\bf c}}  \pi_{\bf c}\widehat{g}\,)\big\|_{L^q(\mathbb R^3 \times I)} 
\le C_\delta \lambda^{C-\frac 32\nu D}\|g\|_p, \quad 1\le p\le q\le \infty.\]

Now it remains to show \eqref{Ie}. We recall $a_e(s,\xi)=\psi_{\!J}(s+\sigma(\xi)) ( 1-\beta_0(2^{-4}\lambda^{\frac13-\nu}s))$.  
Since  $|s| \ge 2^4 \lambda^{\nu-\frac13}$ on 
$\supp a_e(\cdot, \xi)$ and   $|R_\gamma(\xi)|\le 2\lambda^{2\nu-\frac 23}$ for $\xi\in \supp \pi_{\bf c}$, 
by \eqref{taylor-phi} and \eqref{taylor-phi2} it follows that  
$ C_1   \lambda |s|^{2} \le   |\partial_s \phi(s,\xi)| \le C_2  \lambda |s|^{2}$
 and 
\begin{align}
\label{derivative-bound}
\begin{aligned}
C_3\lambda |s|^{3-k} \le  |\partial_s^k \phi(s,&\xi)|\le   C_4 \lambda |s|^{3-k}, \quad \ \ k=2,3,  
\\
| \partial_s^k \phi(s,\xi)|&\le C_5 \vep_\circ \lambda, \quad\quad  \qquad  \ 4\le k\le D
\end{aligned}
\end{align} 
for some positive constants $C_1,\dots, C_5$. Thus,  
  noting  $|\partial_s^k a_e(s,\xi) |  \le C_\delta \lambda^{(\frac 13-\nu)k}$  for $0\le k\le D$, 
    we have  
\Be 
 \label{int-by}
 b_{\ell+1} :=  \frac{|\partial_s^{\ell+1} \phi(s,\xi)|}{|\partial_s \phi(s,\xi)|^{\ell+1}}\le C_\delta \lambda^{-\frac32\nu (\ell+1)},  
   \quad \ 
  b'_{\ell}:=\frac{|\partial_s^{\ell} a_e(s,\xi)|}{|\partial_s \phi(s,\xi)|^{\ell}}\le C_\delta \lambda^{-3\nu \ell}  
    \Ee 
   for $\ell \ge 1$ 
 if $\xi\in \supp \pi_{\bf c}$ and $|s| \ge 2^4 \lambda^{\nu-\frac13}$.   
 After integration  by parts  $D-1$ times  we see that 
$ |\mathcal I_e(t,\xi)|$  is bounded  by  a finite sum of  the integrals $ C\int  
\prod_{j=1}^{m}  \mathcal M_{\ell_j} ds$ where $\mathcal M_{\ell}\in \{b_{\ell+1}, b'_{\ell}\}$, $\sum_{j=1}^m\ell_j=D-1$, and $\ell_j\ge 1$. 
Using \eqref{int-by}  we  get 
$ |\mathcal I_e(t,\xi)|\le C_\delta \lambda^{-\frac 32\nu D}$ 
for $\xi\in \supp \pi_{\bf c}$.  Furthermore, since $\partial_s^k\partial_\xi^\alpha \phi$, $\alpha\neq 0$ are bounded, 
 the same argument  shows   \eqref{Ie}. 
\end{proof}

\subsubsection*{Multipliers of  $ A^\gamma[\psi_{\!J}] \mathcal P_{\!\mathbf {o}}^1$ and $A^\gamma[\psi_{\!J}] \mathcal P_{\!\mathbf {o}}^0$}  
We  obtain  similar expansions  for $m_{\!J}\pi_{\bf o}^j$, $j=0,1$.  As we shall see, 
$m_{\!J}\pi_{\bf o}^0$ is decaying  rapidly as $\lambda\to\infty$ (see \eqref{m-near} below).   
We    concentrate on the case  $\xi\in \supp  \pi_{\mathbf o}^1$ for the moment.  

Let $\rho_1\in \mathrm C_0^\infty([2^{-5}, 2^5])$, $\rho_0\in \mathrm C_c^\infty([0, 2^{-4}))$, and
$\rho_2\in \mathrm C^\infty((2^{4}, \infty))$   such that
$\rho_1=1$ on $[2^{-4}, 2^4]$ and   $\rho_0+\rho_1+\rho_2=1$ on $[0,\infty)$.  For $j=0,1,2,$ we set 
\begin{align*}
  a_j(s,\xi) &= \psi_{\!J} (s+\sigma(\xi)) \rho_j\big(R^{-1/2}_\gamma(\xi)|s|\big),
  \\
    \mathcal I_j(t,\xi) &= (2\pi)^{-3}\int e^{-it\phi(s,\xi)} a_j(s,\xi)\,ds, 
  \end{align*}
and then we have 
\Be 
\label{decomp-m} m_{\!J}(t,\xi)   
=   e^{-it \Phi^{\bf c}(\xi)}\big( \mathcal I_0(t,\xi)  +\mathcal I_1(t,\xi)+\mathcal I_2(t,\xi)\big).
 \Ee
 The main term is $\mathcal I_1$ while $ \mathcal I_0$ and $\mathcal I_2$ are rapidly decaying  as $\lambda\to \infty$ (see \eqref{b0b2-near} below).
 The second derivative of the phase function does not vanish on $\supp a_1(\cdot, \xi)$, so 
 we may apply the method of stationary phase  for   $\mathcal I_1(t,\xi)$. 
 For the purpose  we set 
\Be 
\label{tilde-phi}\widetilde \phi(s,\xi)
             =  L^{-1}(\xi)  \phi\big(R^{1/2}_\gamma(\xi) s,\xi\big)
             \Ee
where $L(\xi)=\Lambda_\gamma(\xi)R_\gamma(\xi)^{\frac 32}$,  and set
\begin{align*}
 a^\pm(s, \xi)
              &=\psi_{\!J}\big(R^{1/2}_\gamma(\xi)s+\sigma(\xi)\big) \rho_1(\pm s),
              \\
\mathcal I_1^\pm (t,\xi)&= (2\pi)^{-3} R^{1/2}_\gamma(\xi)
\int e^{-it L(\xi)  \widetilde \phi(s,\xi)}   a^\pm(s, \xi)\,ds.           
\end{align*}
By scaling $s \rightarrow R^{1/2}_\gamma(\xi) s$ we  have
\begin{equation}\label{m12} 
\mathcal I_1(t,\xi)=\mathcal I_1^+(t,\xi)+\mathcal I_1^-(t,\xi). 
\end{equation}

We try to find the stationary  points of  the function  $\widetilde \phi(\cdot, \xi)$, which give rise to two different phase functions $\Phi^\pm$ (see \eqref{phase} below). 
As we shall  see later, it is important for application of the multilinear restriction estimate how smooth these phase functions  are. So, we treat the matter carefully.

\begin{lem} 
\label{taupm} There are $\tau_+, \tau_- \in   \mathrm C^{D-4}(\mathbb A^\ast_\lambda\times [-\delta^{10}, \delta^{10}] )$, homogeneous of degree zero, such that  
$ \pm \tau_\pm(\xi,\theta)\in  [2^{-1}, 2]$ and, if $R_\gamma(\xi)\ge 0$, 
\Be\label{t0}
 \partial_s \widetilde \phi\big(  \tau_\pm(\xi, R^{1/2}_\gamma(\xi)), \xi \big)=0.
 \Ee
\end{lem} 

\begin{proof} Recalling \eqref{taylor-phi}, we  set  \[\Theta_0(s,\xi)=s^{-3} \Theta(s,\xi),\] 
which is homogeneous of degree zero in $\xi$. One can see  $\Theta_0 \in \mathrm C^{D-3}([-1/2,1/2]\times \mathbb A^\ast_\lambda)$ without difficulty  because 
$
\Theta_0(s,\xi)= (s/{3!}) \int_0^1 {(1-t)^3 \gamma^{(4)}( st+\sigma(\xi)) \cdot   \xi} {\, \Lambda_\gamma^{-1}(\xi)}  dt 
$ by Taylor's theorem with integral remainder. 
Then we  consider the function 
 \[  \widetilde \phi_0(s,\xi,\theta)=- s+ \frac{s^3}{3!} 
+  s^3 \Theta_0(\theta s,\xi)\] 
with $(s,\xi,\theta)\in  \Omega^\pm:= ( \pm [2^{-5}, 2^5]) \times \mathbb A^\ast_\lambda\times [-\delta^{10}, \delta^{10}]$.  
It is clear that $\widetilde \phi_0\in \mathrm C^{D-3}(\Omega^\pm)$.

Since $\Theta_0$, $\partial_s \Theta_0$ and $\partial_s^2 \Theta_0$ are $\mathcal O(\vepc)$  as can be seen  using \eqref{taylor-phi} and \eqref{taylor-phi2}, 
we have
$\partial_s \widetilde \phi_0(s,\xi,\theta)=-1+  {s^2}/2+ \mathcal O(\vepc)$ and $\partial_s^2 \widetilde \phi_0(s,\xi,\theta)= s+ \mathcal O(\vepc) $. We now note that  $\partial_s \widetilde \phi_0(\cdot, \xi, \theta)$ has two distinct zeros which are respectively close to $\sqrt 2$ and $-\sqrt 2$, thus by the implicit function theorem  there are 
$\tau_+(\xi,\theta)$ and $\tau_-( \xi, \theta)$  such that   $\partial_s \widetilde \phi_0\big(  \tau_\pm(\xi, \theta), \xi,\theta \big)=0$ 
and  $ \pm \tau_\pm(\xi,\theta)\in  [2^{-1}, 2]$ if $\vepc$ is small enough.  Additionally,  $\tau_+$ and $\tau_-$ are $D-4$ times continuously differentiable  since so is $\partial_s \widetilde \phi_0$. 
By \eqref{taylor-phi} and \eqref{tilde-phi} we have     $\widetilde \phi_0( s,\xi,R^{1/2}_\gamma(\xi))=  \widetilde \phi(s,\xi)$, thus it follows that  $\partial_s  \widetilde \phi_0( s,\xi,R^{1/2}_\gamma(\xi))= \partial_s  \widetilde \phi(s,\xi)$ when $R_\gamma (\xi)\ge 0$.   
Therefore  we obtain  \eqref{t0}.
\end{proof}

We  set 
\[s_\pm(\xi)=R^{1/2}_\gamma(\xi)\tau_\pm \big(  \xi,R^{1/2}_\gamma(\xi)\big).\]
Then from \eqref{tilde-phi} 
  it follows  that  $\displaystyle{\gamma'\big(s_\pm(\xi)+\sigma(\xi)\big)\cdot \xi =0}$. We  define
\begin{align}\label{phase}
\Phi^{\pm}(\xi)=  \gamma\big(s_\pm(\xi)+\sigma(\xi)\big)\cdot \xi
\end{align}
for $\xi\in  \mathbb A^\ast_\lambda\cap \{ \xi: R_\gamma(\xi)\ge 0\}$. If $R_\gamma(\xi)=0$ for some $\xi$,   $\nabla\Phi^{\pm}(\xi)$ may not exist because $R^{1/2}_\gamma$ is not differentiable at $\xi$. However, 
$\nabla\Phi^{\pm}$ can be defined to be a continuous function on  $ \mathbb A^\ast_\lambda\cap \{ \xi: R_\gamma(\xi)\ge 0\}$. Indeed, 
differentiating  \eqref{phase} gives 
\Be
\label{gradient}
 \nabla\Phi^{\pm}(\xi)= \gamma\big( s_\pm(\xi)+\sigma(\xi)\big) 
 \Ee
  if  $R_\gamma(\xi)>0.$  Thus   $\nabla\Phi^{\pm}$ becomes continuous on  $ \mathbb A^\ast_\lambda\cap \{ \xi: R_\gamma(\xi)\ge 0\}$  if we set $ \nabla   \Phi^{\pm}(\xi)= \gamma\big(\sigma(\xi)\big)$ when $R_\gamma(\xi)=0$ since $\gamma$, $\sigma$ are continuous.

We  define  the adjoint restriction operators $\mathcal T^\pm_\lambda$ by 
\begin{align*}
\mathcal T^{\pm}_\lambda g(x,t)
=
\int_{\mathcal C_\lambda^{\bf o}(\delta)}  e^{i(x\cdot \xi-t\Phi^{\pm}(\xi))} g(\xi)\,d\xi,
\end{align*} 
where $
\mathcal C_\lambda^{\bf o} (\delta):= \{\xi\in  \mathbb A^\ast_\lambda: 0\le R_\gamma(\xi) \le 2\delta^{100} \}$. 
 Putting together the  discussion so far with the method of stationary phase  we can obtain

\begin{lem}\label{op3} Let $0<\nu \ll1$, $M=[\frac {D-1}3]$, and $J\in \mathfrak J_\circ(\delta)$.   
 Suppose  $\gamma \in \mathfrak C^D(\varepsilon_\circ)$, $\psi_{\!J}\in \mathfrak N^D(J)$, and $\widehat f$ is supported on $\mathbb A_\lambda$.
Then, we have 
\Be
\label{eq:whole}
 A^\gamma[\psi_{\!J}] (\mathcal P_{\!\mathbf {o}}^1+ \mathcal P_{\!\mathbf {o}}^0) f
=\sum_{\pm}\sum_{\ell=0}^{M-1} t^{-\frac{2\ell+1} 2}
\mathcal T^{\pm}_\lambda \big( \gamma_\ell^\pm \pi_{\bf o}^1\widehat{f}\,\,\big) 
+
\mathcal E_{\mathbf o}  f, \quad  t\in I,
\Ee
and the following hold with $C$ and  $C_\delta$  independent of   $\gamma$ and $\psi_{\!J}$:  
\begin{align}\label{dec-b1}
|\gamma_{\ell}^{\pm}(\xi)|\le C_\delta \lambda^{-\frac 13-\frac \nu2}
\lambda^{-3 \ell \nu}
\end{align}
 for $0\le \ell\le M-1$  and 
\Be
\label{error-est} 
\| \mathcal E_{\mathbf o} f\|_{L^q(\mathbb R^3 \times I)} \le C_\delta \lambda^{C- 3\nu M} \|  f\|_p, \quad 1 \le p\le q \le \infty.
\Ee
\end{lem}

It should be noted that the expansion in \eqref{eq:whole}  is obtained only on the support of $ \pi_{\bf o}^1$ but not on  the larger set
$\mathcal C_\lambda^{\bf o}(\delta)$.

We now proceed to apply to $\mathcal I_1^\pm $ the  method of stationary phase.  We first note that $\supp a^\pm(\cdot, \xi) \subset \pm[2^{-5}, 2^5]$ and, as seen in the above,  the phase $ \widetilde \phi(\cdot,\xi)$ has 
the stationary points $\tau_\pm(  \xi, R^{1/2}_\gamma(\xi))$ while 
$\partial_s^2\widetilde \phi(\cdot,\xi)=s+\os(\vepc)$ for $\xi\in  \mathbb A^\ast_\lambda\cap \{ \xi: 0\le R_\gamma(\xi)\le 2\delta^{100}\}$.    
We also note that   $|L(\xi)| \ge 2^{-1} \lambda^{3\nu}$ for $\xi\in \supp \pi_{\bf o}^1$ and  that  $ L(\xi) \widetilde \phi(  \tau_\pm(  \xi, R^{1/2}_\gamma(\xi)) ,\xi)
             =   \gamma\big( s_\pm(\xi)  +\sigma(\xi)\big)\cdot \xi -\Phi^{\bf c}(\xi)$. 
             Bring all these observations together, we now apply \cite[Theorem 7.7.5]{H} (also see  \cite[Theorem 7.7.6]{H})  and obtain 
\begin{align}
\label{asymp}
{\mathcal I_1^\pm}(t,\xi)
=e^{it (\Phi^{\bf c}(\xi) -\Phi^{\pm}(\xi) )} R^{1/2}_\gamma(\xi)
\sum_{\ell=0}^{M-1} d_{\ell}^{\pm}(\xi) & (tL(\xi))^{-\frac12-\ell}
+  e_{M}^\pm (t,\xi) 
\end{align}
for $\xi\in \supp \pi_{\bf o}^1$  where $M=[\frac {D-1}3]$ and  $e_{M}^\pm (t,\xi)=\mathcal O\big( |tL(\xi)|^{-M} \big)$.  
The functions $d_{\ell}^{\pm}$ are  bounded   on the support of  $\pi_{\bf o}^1$ since  so are $\partial_s^k \widetilde \phi$ and $\partial_s^k a^\pm$.

\begin{proof}[Proof of Lemma \ref{op3}] 
Recalling \eqref{decomp-m} and  \eqref{m12},  we   write 
       \[ m_{\!J} (\pi_{\bf o}^1+\pi_{\bf o}^0)= e^{-it \Phi^{\bf c}} (\mathcal I_1^+ + \mathcal I_1^{-}) \pi_{\bf o}^1+  e^{-it \Phi^{\bf c}} (\mathcal I_0+\mathcal I_2)\pi_{\bf o}^1+ m_{\!J}  \pi_{\bf o}^0.  \]    
 Using  \eqref{asymp},  we now put 
 \begin{align*}
  \mathcal E(t,\cdot)&= e^{-it \Phi^{\bf c}}  \big(   e_{M}^+ (t,\cdot) + e_{M}^- (t,\cdot) \big)  \pi_{\bf o}^1
+   e^{-it \Phi^{\bf c}} \big(\mathcal I_0(t,\cdot)+\mathcal I_2(t,\cdot)\big)\pi_{\bf o}^1+ m_{\!J}(t,\cdot)  \pi_{\bf o}^0,
  \end{align*}
and then we set $\mathcal E_{\bf o} f= (2\pi)^{3} \mathcal F^{-1}_\xi(  \mathcal E(t,\cdot)\widehat f\,)$ and $\gamma^\pm_\ell(\xi)=R^{1/2}_\gamma(\xi)d_{\ell}^{\pm}(\xi)(L(\xi))^{-\frac12-\ell}.$ 
Thus we have \eqref{eq:whole}.  Recalling $L(\xi)=\Lambda_\gamma(\xi)R_\gamma(\xi)^{\frac 32}$, we  see  \eqref{dec-b1} holds since  $C\lambda^{2\nu-2/3}\le  |R_\gamma(\xi)|$ and $d_\ell^\pm$ are bounded on  $\supp \pi_{\bf o}^1$. 

To show \eqref{error-est}  we  use
   the following: 
\Be
\label{m-near}
 |\partial_\xi^\alpha m_{\!J}(t,\xi)|\le C_\delta \lambda^{-\frac 32\nu (D-|\alpha|)}, \qquad \xi\in \supp  \pi_{\mathbf o}^0, 
 \Ee
and  
\Be 
\label{b0b2-near} 
\begin{aligned}
 |\partial_\xi^\alpha \mathcal I_0(t,\xi)|&\le  C_\delta \lambda^{-\frac 32\nu (D-|\alpha|)},
 \\
   |\partial_\xi^\alpha \mathcal I_2(t,\xi)| &\le  C_\delta \lambda^{-\frac 32\nu (D-|\alpha|)},
   \end{aligned}
   \qquad   \xi\in \supp  \pi_{\mathbf o}^1.
\Ee 
 Assuming this for the moment we obtain \eqref{error-est}.  
 Note that  $|\partial_\xi^\alpha \Phi^{\bf c}|\le C\lambda^{1-|\alpha|}$ and $|\partial_\xi^\alpha e_M^\pm |\le C_1\lambda^{C-3\nu M}$ for $|\alpha|\le 4$. Combining this,  \eqref{m-near} and \eqref{b0b2-near}  for $|\alpha|\le 4$  and using the  estimate \eqref{easy1},  we get  \eqref{error-est}  in the same manner  as before.

To complete the proof, we are left to prove  \eqref{m-near} and  \eqref{b0b2-near}. Let us first consider \eqref{m-near}, which is easier. Since  $R_\gamma(\xi)\le -\lambda^{2\nu-\frac23}$ for $\xi\in \supp  \pi_{\mathbf o}^0$,  by  \eqref{taylor-phi}  we see that 
$
|\partial_s \phi| \ge C_1 \lambda \big( -R_\gamma(\xi) + s^2(1/2-\varepsilon_\circ |s|)  \big)
\ge C_2 \lambda\max(s^2, \lambda^{2\nu-\frac23})$  for some $C_1, C_2>0$. Combining this with \eqref{derivative-bound}, 
we have  \eqref{int-by} for $\ell\ge 1$ with $a_e(s,\xi)$ replaced by $\psi_J(s+\sigma(\xi))$.  Thus integration by parts gives $|m_{\!J}(t,\xi)|\le C_\delta\lambda^{-\frac 32\nu D}$ since $R_\gamma(\xi)\le -\lambda^{2\nu-\frac23}$.  
The same argument  also works for $\partial_\xi^\alpha m_{\!J}(t,\xi)$,  so we obtain 
\eqref{m-near}.

We  now show \eqref{b0b2-near} only with $\alpha=0$, and 
 the derivatives $\partial_\xi^\alpha \mathcal I_0$ and $\partial_\xi^\alpha \mathcal I_2$ can be handled likewise. 
   We consider  $\mathcal I_0$ first. By  \eqref{taylor-phi} we have $ |\partial_s \phi|\ge C \lambda R_\gamma(\xi)$ for $|s| \le 2^{-4}R_\gamma^{1/2}(\xi)$. 
Combining this with \eqref{derivative-bound}, we get the first estimate in \eqref{int-by}
for $\ell \ge 1$  when $|s| \le 2^{-4}R_\gamma^{1/2}(\xi)$ because $\lambda^{2\nu-\frac 23} \le R_\gamma(\xi) $ for $\xi\in \supp  \pi_{\mathbf o}^1$. 
Note that $|\partial_s^\ell a_0(s,\xi)|\le C_\delta R_\gamma^{-\ell /2}(\xi)$, hence for $\ell\ge 1$ we have 
the second estimate in \eqref{int-by}  with $a_e$ replaced by $a_0$.
Therefore repeated integration by parts gives the  estimate for $\mathcal I_0$.  We can handle  $\mathcal I_2$ in the same manner. 
 Since $|s| \ge 2^{4}R_\gamma^{1/2}(\xi)$, by \eqref{taylor-phi}  we have  $C_1   \lambda |s|^{2}\le |\partial_s \phi(s,\xi)|\le C_2   \lambda |s|^{2}$ and obviously $|\partial_s^\ell a_2(s,\xi)|\le C_\delta R_\gamma^{- \ell/2}(\xi)$.  So,   we get the estimate \eqref{int-by} for $|s| \ge 2^{4}R_\gamma^{1/2}(\xi)$ and $\ell\ge 1$ while $a_e$ is replaced by $a_2$. Thus  integration by parts gives  the estimate for $\mathcal I_2$.   
\end{proof}

In contrast to $\Phi^{\bf c}$  the 2nd derivatives of $\Phi^{\pm}$ are no longer bounded. However, a computation with $\gamma=\gamma_\circ$ \footnote{
If $\gamma=\gamma_\circ$,  $
\Phi^{\bf c}(\xi)=-{\xi_1\xi_2}/{\xi_3}+{\xi_2^3}/(3\xi_3^2) $ and $
\Phi^\pm(\xi)= \Phi^{\bf c}(\xi) 
\mp 3^{-1}  {\xi_3} \big( {\xi_2^2}/{\xi_3^2} - 2{\xi_1}/{\xi_3} \big)^{3/2}.
$
}  leads us to expect  that  $\Phi^{\pm}\in \mathrm C^{1,1/2}$.  What follows  shows   this holds true for $\gamma\in \classe$.

\begin{lem} 
\label{lem:holder} For $ \xi_1, \xi_2\in \mathcal C_1^{\bf o} (\delta)$, there is a constant $C$ independent of $\gamma$  such that 
\Be
\label{holder} 
 |\nabla  \Phi^{\pm}(\xi_1) -\nabla \Phi^{\pm}(\xi_2)|\le C| \xi_1-\xi_2|^\frac12.
 \Ee
\end{lem}

\begin{proof} 
Let us set $\tau_0^\pm (\xi) =\tau_\pm ( \xi, R^{1/2}_\gamma(\xi))$, so  $s_\pm(\xi) =R^{1/2}_\gamma(\xi) \tau_0^\pm (\xi) $. 
Using  \eqref{gradient} and applying the mean value inequality,   one can easily  see
\[ |\nabla  \Phi^{\pm}(\xi_1) -\nabla \Phi^{\pm}(\xi_2)|\le C  |  s_\pm(\xi_1)- s_\pm(\xi_2)| + C|\sigma(\xi_1) -\sigma(\xi_2) |.\] 
Since $\sigma \in \mathrm C^{D-2}(\mathbb A_\lambda^*)$ from Lemma \ref{sigma}, we only  have to consider the first one   on the right hand side, which is in turn  bounded by 
\[ |  R^{1/2}_\gamma(\xi_1)- R^{1/2}_\gamma(\xi_2)||\tau_0^\pm(\xi_1)|+  R^{1/2}_\gamma(\xi_2)| \tau_0^\pm(\xi_1)- \tau_0^\pm(\xi_2)|. \]
It is easy to see  $|  R^{1/2}_\gamma(\xi_1)- R^{1/2}_\gamma(\xi_2)|\le C |\xi_1-\xi_2|^\frac12$. 
Since $ \tau_\pm$ is $D-4$ times continuously differentiable in a region containing $\mathcal C_1^{\bf o} (\delta)$ (Lemma \ref{taupm})  and 
$ \tau_0^\pm(\xi)=  \tau_\pm( \xi, R^{1/2}_\gamma(\xi))$, by the mean value inequality it follows that  
$| \tau_0^\pm(\xi_1)- \tau_0^\pm(\xi_2)|\le C |  R^{1/2}_\gamma(\xi_1)- R^{1/2}_\gamma(\xi_2)|+  C |\xi_1-\xi_2|. $ Consequently, we get the inequality \eqref{holder}. 
\end{proof}

\subsection{Multilinear restriction estimate}
In this section, we obtain a form of multilinear restriction estimate, which we need to prove \eqref{global}.  The surfaces associated with $\Phi^{\bf c}$ and $\Phi^{\pm}$  have some curvature property, so it is possible to get an $L^2$-$L^q$ smoothing estimate using the typical $TT^\ast$ argument. However, the consequent estimate is not so strong enough to be useful for controlling the maximal operator. Instead, we utilize $4$-linear estimates which we deduce from the multilinear restriction estimate under transversality assumption (\cite{BCT}).

\subsubsection*{Multilinear restriction estimate for $\mathrm C^{1,\alpha}$ hypersurfaces}   For the adjoint restriction estimate,  the surfaces  are typically assumed to be compact and twice continuously  differentiable. The same assumption was also made for the multilinear restriction estimate in \cite[Theorem 1.16]{BCT}, but  the phase   functions $\Phi^\pm$ 
 no longer have bounded second derivatives. Nevertheless, it is not difficult to see that the argument in \cite{BCT}  continues  to work with $\mathrm C^{1,\alpha}$ surface, $\alpha>0$.  
 It seems to the authors that there is no proper reference concerning this matter, so we provide a brief discussion on the multilinear restriction estimate for the less regular $\rm C^{1,\alpha}$ surfaces.

For  $k=1,\dots d$, let  $U_k$ be a compact subset of an open set   $U_k' \subset \mathbb B^{d-1}(0, 2^2)$
and $\Phi_k$ be a real valued function  on 
$U_k'$ which satisfies 
$
\| \Phi_k \|_{\mathrm C^{1,\alpha}(U_k')} \le B
$ 
for some $0<\alpha \le 1$. Let us set 
\[
T_kg_k(x,t)= \int_{U_k} e^{i(x\cdot \xi -t\Phi_k(\xi))} g_k(\xi)\,d\xi.
\]

\begin{thm}
\label{holder-rest}
Let  $d \ge2$, $\theta\in (0, 1]$,  and let  $N_k(\xi)=\frac{( \nabla \Phi_k(\xi),1)}{|(\nabla \Phi_k(\xi),1)|} $.  
Suppose $|\det( N_1(\xi_1),\dots, N_d(\xi_d))| \ge \theta $
for $\xi_k \in U_k$, $k=1,\dots,d$.  Then, for any $\vep>0$, 
\begin{align}\label{multi-tg}
\big\| \prod_{k=1}^d T_kg_k \big\|_{L^{\frac 2{d-1}}( \mathbb B^d(0,R))}
\le C_\vep(\theta) R^\varepsilon
\prod_{k=1}^d \| g_k \|_{L^2(U_k)}
\end{align}
holds whenever $R \ge 1$.  The constant $C_\varepsilon(\theta)$ takes the form $C\theta^{-C_\vep}$ for some constants $C$,  $C_\vep>0$.
\end{thm}

When $d=2$, the theorem holds true with a $\mathrm C^1$ curve even with a Lipschitz curve, but it is unknown whether the same continues to be true in higher dimensions. 
Once one makes  a couple of crucial  observations concerning  $\mathrm C^{1,\alpha}$ surfaces,  
it is not difficult to prove Theorem \ref{holder-rest} through routine adaptation of the arguments in \cite[Proposition 2.1]{BCT}. Instead of reproducing them in detail,  we provide 
a sketch of the proof. We refer the reader to \cite{BCT, Bennett} for the details.  

For the proof of  Theorem \ref{holder-rest}, first of all,  we observe that 
\Be
\label{holder-con}
 | \Phi_k(\xi +h) -\Phi_k(\xi)-\nabla  \Phi_k(\xi)\cdot h|\le C B|h|^{\alpha+1}
 \Ee
for $\xi+h, \xi\in U_k$.  
If $\Phi_k$ is assumed to be in $\mathrm C^{1,\alpha}(U_k)$ instead of $\mathrm C^{1,\alpha}(U_k')$,  this can  not  be completely clear. In such a case we need to impose an additional  condition such that $U_k$ has  a $\mathrm C^{1,\alpha}$ boundary (e.g. see \cite[pp. 136--137]{GT}).  On the other hand, if  $U_k$ is convex, \eqref{holder-con} is a simple consequence of the mean value theorem. 
Since $U_k$ is compact,  there is a positive  number $\rho_k$ such that  
$x,y$ are contained in a ball  which is a subset of $U'_k$  whenever $x,y\in U_k$ and  $|x-y|\le \rho_k$.  Therefore, we get \eqref{holder-con} for $|h|\le \rho_k$ and this is enough to show \eqref{holder-con} for any $\xi, \xi+h\in U_k$ because $U_k$ is compact and $\nabla \Phi_k$ is continuous.

\subsubsection*{Sketch of proof} 
Let us denote $\Sigma_k=\{ (\xi, -\Phi_k(\xi)): \xi\in U_k\}$. We consider the estimate 
\begin{align}\label{G2}
\big\| \prod_{k=1}^d \widehat{G_k} \big\|_{L^{\frac2{d-1}}(\mathbb B^d(0,R))}
\le
C_0 R^{-\frac d2}
\prod_{k=1}^d \|G_k\|_{L^2(\mathbb R^d)}
\end{align}
for  $R\ge 1$ when $G_k$ is supported in  $\Sigma_k(1/R):=\{ (\xi,\tau)\in \mathbb R^{d-1} \times \mathbb R: \dist((\xi, \tau), \Sigma_k)<1/R \}$. The estimate \eqref{multi-tg} is equivalent to \eqref{G2} with
$C_0= CR^\varepsilon$ (see \cite{BCT}).  Let $\mathcal C(R)$ be the infimum of $C_0$ with which \eqref{G2} holds.  
The key part of the proof  is to establish the implication 
\Be \label{imply} \mathcal C(R)\le R^b  \implies   \mathcal C(R)\le C(\theta,\vep) R^{\frac b{1+\alpha}+\vep} \Ee
for any $\vep>0$ where $b$ is a positive constant.   Via iteration, the exponent of $R$ can be suppressed to be arbitrarily  small  and 
hence we get the estimate \eqref{multi-tg}.

 Using \eqref{holder-con}  we see that the set  $\Sigma_k(1/R)\cap \mathbb B^d(\zeta, R^{-1/(1+\alpha)})$, $\zeta\in  \Sigma_k$ is contained in a $C/R$ neighborhood of 
the tangent plane to  $\Sigma_k$ at $\zeta$.  
Thus $\Sigma_k(1/R)$ can be covered with a collection $\{\mathfrak R_{j}^k\}$ of  
finitely overlapping rectangles    of dimensions  about  $R^{-1}\times R^{-1/(1+\alpha)}\times \cdots\times R^{-1/(1+\alpha)}$ which are essentially  tangential to $\Sigma_k(1/R)$. These rectangles provide a decomposition of  
$G_k=\sum_{j} G_{j}^k$ while $\supp G_{j}^k\subset \mathfrak R_{j}^k$.  Thus, after applying the assumption $\mathcal C(R)\le R^b$ to  the integrals over the balls  of radius $R^{1/(1+\alpha)}$ which finitely overlap and cover $\mathbb B^d(0,R)$,  one can get the implication \eqref{imply} using the multilinear Kakeya estimate \cite{BCT, Guth} for the transversal collection $\mathfrak T_k$ of  the tubes of  width  $R^{1/(1+\alpha)}$ and length $R$ which have their axes parallel to the normal vector of the surface $\Sigma_k $.   \qed

\smallskip

Making use  of Theorem \ref{holder-rest} we obtain the following.

\begin{prop}\label{prop:multi0}
Let  $\theta_1, \dots, \theta_4\in\{\bf c, +,-\}$ and let $J_k \in \mathfrak J_\circ(\delta)$,  $1\le k\le 4$.   Suppose that  $\gamma \in \mathfrak C^D(\varepsilon_\circ)$ and  $\dist(J_\ell, J_k)\ge \delta$, $\ell\neq k$. 
Then, for $\varepsilon>0$ and $R\ge 1$,   there is a constant  $C_{\varepsilon}$ such that 
\begin{align}\label{T4-}
\Big\| \prod_{k=1}^4
|\mathcal T_{1}^{\theta_k} \big( {\widetilde \chi}_{\mathcal R_{\!J_k}}(\lambda\cdot) g_k\big) |^{\frac 14} 
\Big\|_{L^{\frac 83}(\mathbb B^4(0,R))}
&\le C\delta^{-C_\vep} R^\varepsilon \prod_{k=1}^4  
\| g_k\|_2^{\frac 14}.
\end{align}
\end{prop}

\begin{proof}
We begin with recalling that ${\widetilde \chi}_{\mathcal R_{J_k}}(\lambda \cdot)$ is supported in $\lambda^{-1}\mathcal R_{J_k}(2^6) $ and  that $|R_\gamma(\xi)|\le 2\delta^{100}$ if $\xi\in \mathcal C_1^{\bf c} (\delta)$ or $\mathcal C_1^{\bf o} (\delta)$.  Since $\nabla_\xi \Phi^{\bf c}(\xi)=\gamma(\sigma(\xi))
+\gamma'(\sigma(\xi))\cdot \xi\, \nabla \sigma(\xi)$,   we have
$
\nabla_\xi \Phi^{\bf c}(\xi)
=\gamma(\sigma(\xi))+\mathcal O(\delta^{100})
$
for   $\xi\in \mathcal C_1^{\bf c} (\delta)$. 
If $\xi\in \mathcal C_1^{\bf o} (\delta)$, by \eqref{gradient}  we have
$
\nabla_\xi \Phi^{\pm}(\xi)
=\gamma(\sigma(\xi))+\os(2^2 \delta^{50})
$
because  $|R_\gamma(\xi)|\le 2\delta^{100}$. 
  Thus  
  \[ \mathrm N_k(\xi):= |(\nabla \Phi^{\theta_k}(\xi), 1)|^{-1} {(\nabla \Phi^{\theta_k}(\xi), 1)}\] which is normal  to the surface $(\xi,- \Phi^{\theta_k}(\xi))$ satisfies 
\begin{align*}
\mathrm N_k(\xi)=\frac{(\gamma(\sigma(\xi)), 1)}{\sqrt{|\gamma(\sigma(\xi))|^2+1}  }  + \os( 2^3
\delta^{50}), \quad  \xi\in \mathcal C_1^{\theta_k}(\delta),  \quad k=1,\dots, 4,
\end{align*}
where we denote  $\mathcal C_1^{\pm}(\delta)=\mathcal C_1^{\bf o}(\delta)$. 

Let $\xi_k\in \lambda^{-1}\mathcal R_{J_k}(2^6)\cap \mathcal C_1^{\theta_k}(\delta)$, $k=1,\dots,4$. 
Then we have $\sigma(\xi_k)\in [-3\cc, 3\cc]$ since $J_k\subset  (1+2\cc) J_\circ$.   
Let $ \Gamma$  denote  the matrix whose $k$-th column is the vector $\big(\gamma(\sigma(\xi_k)), 1\big)$, $k=1,\dots, 4$. 
 By the generalized mean value theorem
(see for example \cite[Part V, Ch.1, 95]{PolyaSzego}) 
there exists $u_k \in [-3\cc, 3\cc]$ such that
\begin{align*}
\det\Gamma =
\det\begin{pmatrix} \gamma(u_1) & \gamma'(u_2 ) & \gamma''(u_3 ) &\gamma'''(u_4 ) \\ 1  & 0   & 0 &  0 \\ \end{pmatrix}
\prod_{1\le \ell<k \le 4} |\sigma(\xi_\ell) -\sigma(\xi_k)|.
\end{align*}
Since $\gamma \in \mathfrak C^D(\varepsilon_\circ)$ and $u_1,\dots, u_4 \in [-3\cc, 3\cc]$,
the determinant on the right hand side has its absolute value $ 1+\os(\vepc)$ regardless of $\gamma$ (for example see \eqref{e123}). 
On the other hand,  using \eqref{SJ2} with $s=c_{J_k}$,  for $\xi_k\in \lambda^{-1}\mathcal R_{J_k}(2^6)\cap \mathcal C_1^{\theta_k}(\delta) $ we have 
$ |c_{J_k}-\sigma(\xi_k) |\le  2^{-2}\delta$ with a small enough $\vepc$,  and  we also have
$|c_{J_\ell}-c_{J_k}|\ge (1+2\cc) \delta $, $\ell\neq k$ because  $\dist(J_\ell, J_k)\ge \delta$. So, $ |\sigma(\xi_\ell)-\sigma(\xi_k) |> 2^{-1} \delta$
if $\ell\neq k$, and   we thus have 
 $\prod_{1\le \ell < k \le 4} |\sigma(\xi_\ell)-\sigma(\xi_k) |> 2^{-6}   \delta^6$. 
Consequently, we obtain 
\[
 |\!\det (\mathrm N_1(\xi_1), \dots, \mathrm N_4(\xi_4))| 
 > 2^{-7} \delta^6
\]
provided that  $\xi_k\in \lambda^{-1}\mathcal R_{J_k}(2^6)\cap \mathcal C_1^{\theta_k}(\delta)$ for $k=1,\dots,4$. That is to say, 
the transversality condition holds uniformly  regardless of the choice of  $\theta_1, \dots, \theta_4\in\{\bf c, +,-\}$.

 We now note that  $\Phi^{\bf c}$ is  continuously differentiable at least twice  in a region containing $\mathcal C_1^{\bf c} (\delta)$ and  that 
 $\| \Phi^{\pm}\|_{\mathrm C^{1,1/2}(\mathcal C_1^{\bf o}(\delta))}\le C$ by Lemma \ref{lem:holder}. 
 To apply Theorem \ref{holder-rest}  we need only to make it sure that $\Phi^{\pm}$ extends as a $\mathrm C^{1,1/2}$ function to an open set containing  $\mathcal C_1^{\bf o}(\delta)$. 
The only part of the boundary which can be problematic is $S:=\{ \xi: R_\gamma(\xi)=0\}\cap \mathcal C_1^{\bf o}(\delta)$ since $\Phi^{\pm }$ is  homogenous and 
$D-4$ times continuously differentiable on 
$\{ \xi: R_\gamma(\xi)=2\delta^{100}\}\cap \mathcal C_1^{\bf o}(\delta)$ (see Lemma \ref{sigma} and  \ref{taupm}).  
 We note that  $R_\gamma(\xi)=0$  if and only if $g(\xi):=\gamma'(\sigma(\xi))\cdot\xi=0$. Since 
   $\nabla g(\xi)=$ $\gamma'(\sigma(\xi))=e_1+\os(6\cc)$  
   for $\xi\in \mathbb A_1^\ast$  by Lemma \ref{sigma} and since 
$g\in \mathrm C^{D-2}(\mathbb A_1^*)$, by the implicit function theorem  it follows that $S$ is a part of  a $\mathrm C^{D-2}$ boundary.   Thus we can extend $\Phi^{\pm }$  to be a $\mathrm C^{1,1/2}$ function across $S$  (e.g., \cite[pp. 136--137]{GT}). Therefore 
we may apply Theorem \ref{holder-rest}  and get the estimate \eqref{T4-}.  \end{proof}

As $\Phi^{\bf c}$, $\Phi^\pm$ are homogeneous of degree $1$,
the following is an immediate consequence of  Proposition  \ref{prop:multi0} by means of scaling and Plancherel's theorem.
\begin{cor}
Under the same assumption as in Proposition \ref{prop:multi0}, for  $\varepsilon>0$, there is a $C_{\varepsilon}=C_\varepsilon(\delta)>0$ such that 
\begin{align}\label{T4}
\Big\| \prod_{k=1}^4
|\mathcal T^{\theta_k}_\lambda \big( {\widetilde \chi}_{\mathcal R_{J_k}} \widehat f_k\big) |^{\frac 14} 
\Big\|_{L^{\frac 83}(\mathbb B^4(0,2^3))}
&\le C_{\varepsilon} \lambda^\varepsilon \prod_{k=1}^4 
\| f_k\|_2^{\frac 14}.
\end{align}
\end{cor}

\subsection{Multilinear estimate for $A^\gamma [\psi_{\!J_k}]  {\mathcal P_{\!\mathbf {n}}}  P_{J_k}$}
We are ready to prove Proposition \ref{global}.  
We first show  quadrilinear estimates 
without weight, from which  we deduce the weighted estimates.

\begin{prop}\label{prop:global}  
Let  $J_k \in \mathfrak J_\circ(\delta)$,  $1\le k\le 4$.   
Suppose that  $\dist(J_\ell, J_k)\ge \delta$, $\ell\neq k$.
If $1/q=5/(8p)+1/16$ and $2\le p \le 6$, then for $\vep>0$,   there  are constants 
$C_\varepsilon = C_{\varepsilon}(\delta)$ and $D=D(\vep)$ such that
\begin{align}\label{A:glob}
\Big\| \prod_{k=1}^4 
|A^\gamma[\psi_{\!J_k}]  {\mathcal P_{\!\mathbf {n}}}  P_{\!J_k}  f_k|^{\frac 14} 
\Big\|_{L^q(\mathbb R^3\times I)}
&\le C_\varepsilon
\lambda^{-\frac 1{3p}-\frac 1{6}+\varepsilon}
\prod_{k=1}^4 \|  f_k\|_{L^p(\mathbb R^3)}^{\frac 14}
\end{align}
holds whenever $\gamma \in \mathfrak C^D(\varepsilon_\circ)$, $\psi_{\!J_k}\in \mathfrak N^D(J_k)$, and  $\widehat f_k$ is supported on $\mathbb A_\lambda$. 
\end{prop}

By the localization argument it is sufficient for the estimate \eqref{A:glob} to show  its local counterpart. 
In fact,  we have 

\begin{lem} 
\label{localtoglobal} Let $1\le p\le q\le \infty$ and $b\in \mathbb R$, and let  $I'\subset I$ be an interval. 
Let $\gamma \in \mathfrak C^D(\varepsilon_\circ)$, $\omega\in \Omega^\alpha$, $0<\alpha\le 4$, and   $\psi_{\!J_k}\in  \mathfrak N^D(J_k)$, $J_k \in \mathfrak J_\circ(\delta)$,  $1\le k\le 4$.  If 
\begin{align}\label{localbd}
\Big\| \prod_{k=1}^4 
|A^\gamma[\psi_{\!J_k}]  {\mathcal P_{\!\mathbf {n}}}  P_{\!J_k}  f_k|^{\frac 14} 
\Big\|_{L^q({\mathbb B^3(0,1)}\times I',\omega)}
&\le B 
\lambda^{b}[\omega]_\alpha^\frac1q
\prod_{k=1}^4 \|  f_k\|_{L^p(\mathbb R^3)}^{\frac 14}
\end{align}
holds for a  large enough $D=D(b)$, then 
we have
\begin{align}\label{globalbd}
\Big\| \prod_{k=1}^4 
|A^\gamma[\psi_{\!J_k}]  {\mathcal P_{\!\mathbf {n}}}  P_{\!J_k}  f_k|^{\frac 14} 
\Big\|_{L^q(\mathbb R^3\times I', \omega)}
&\le C_\delta  B
\lambda^{b} [\omega]_\alpha^\frac1q
\prod_{k=1}^4 \|  f_k\|_{L^p(\mathbb R^3)}^{\frac 14}.
\end{align}
\end{lem}

\begin{proof}  
Let  $K_k(\cdot, t)$  denote  the kernel of the operator  $A^\gamma[\psi_{\!J_k}]  {\mathcal P_{\!\mathbf {n}}}  P_{J_k}$.  We note that the multiplier of  $ {\mathcal P_{\!\mathbf {n}}}  P_{J_k}$ is given by   $m(\xi)= \widetilde \chi_{\mathbb A_\lambda^\ast}(\xi) \beta_0( \delta^{-100} |R_\gamma(\xi)|){\widetilde \chi}_{\mathcal R_{J_k}}(\xi)$ and 
$\|m(\lambda\cdot)\|_{\mathrm C^M}\le C\delta^{-CM}$ for $M\le D-2$.   Since $|\gamma(s)|\le  2(\cc+\vepc)$ for $s\in J_k$,  by Lemma \ref{lem:ker} we have
$ |K_k(x,t)| \le C_\delta   E_M(x)$ for $M\le (D-5)/2$ if $|x|\ge  2$ and $t\in I$. For  $\mathbf  k\in \mathbb Z^3$  set $B_{\mathbf  k}=\mathbb B^3({\mathbf  k,1})$  and $B_{\mathbf  k}'=\mathbb B^3({\mathbf  k},3)$. 
Then  we have 
 \[  |A^\gamma[\psi_{\!J_k}]  {\mathcal P_{\!\mathbf {n}}}  P_{\!J_k}  f|  \le \sum_{\mathbf  k\in \mathbb Z^3}  \chi_{B_{\mathbf  k}}  |A^\gamma[\psi_{\!J_k}]  {\mathcal P_{\!\mathbf {n}}}  P_{\!J_k}  (\chi_{B_{\mathbf  k}'} f) | +  C_\delta E_M \ast | f|.   \]
 Taking $M=4N+9$ above, we combine the inequality with the trivial estimate  $|A^\gamma[\psi_{\!J_k}]  {\mathcal P_{\!\mathbf {n}}}  P_{\!J_k}  g |\le C_\delta \lambda^3 (1+|\cdot|)^{-N} \ast |g|$. Then   we see  
 that $\prod_{k=1}^4 
|A^\gamma[\psi_{\!J_k}]  {\mathcal P_{\!\mathbf {n}}}  P_{\!J_k}  f_k| $ is bounded by 
\[   \!\!  \sum_{\mathbf  k\in \mathbb Z^3}  \chi_{B_{\mathbf  k}}  \prod_{k=1}^4  |A^\gamma[\psi_{\!J_k}]  {\mathcal P_{\!\mathbf {n}}}  P_{\!J_k}  (\chi_{B_{\mathbf  k}'} f_k) |
+ C_\delta \prod_{k=1}^4 (E_N \ast | f_k|).\] 
Since $\|  E_N \ast | f|\|_{L^q(\mathbb R^3\times I', \omega)}\le C   [\omega]_\alpha^{1/q} \lambda^{-N}\|f\|_p$ for $1\le p\le q$, taking a large $N\ge -b,$
we may  disregard the second term. We now use   \eqref{localbd}  to get
\[
\Big\| \prod_{k=1}^4 
|A^\gamma[\psi_{\!J_k}]  {\mathcal P_{\!\mathbf {n}}}  P_{\!J_k} (\chi_{B_{\mathbf  k}'} f_k )|^{\frac 14} 
\Big\|_{L^q(B_\mathbf  k\times I',\omega)}
\le  B
\lambda^{b}[\omega]_\alpha^\frac1q
\prod_{k=1}^4 \|  \chi_{B_{\mathbf  k}'} f_k   \|_{L^p(\mathbb R^3)}^{\frac 14}.
\]
Thus the desired estimate  \eqref{globalbd} follows  by summation over $\mathbf  k$  and H\"older's inequality
since $B_{\mathbf  k}'$  overlap each other at most $6^2$ times. 
\end{proof}

Thanks to  Lemma \ref{localtoglobal},  the proof of Proposition \ref{prop:global} is reduced to showing 
\begin{align}\label{A:glob-2}
\Big\| \prod_{k=1}^4 
|A^\gamma[\psi_{\!J_k}]  {\mathcal P_{\!\mathbf {n}}}  P_{\!J_k}  f_k|^{\frac 14} 
\Big\|_{L^q({\mathbb B^3(0,1)}\times I)}
&\le C_\varepsilon
\lambda^{-\frac 1{3p}-\frac 1{6}+\varepsilon}
\prod_{k=1}^4 \|  f_k\|_{L^p(\mathbb R^3)}^{\frac 14}
\end{align}
for $p,q$ satisfying $1/q=5/(8p)+1/16$ and $2\le p \le 6$. 
Since $\|{\mathcal P_{\!\mathbf {n}}}  P_{\!J_k}  g\|_p\le C_\delta \|g\|_p$ by \eqref{easy2},  using the  estimate \eqref{ps14*} with $p=6$ after H\"older's inequality, we  get 
the estimate \eqref{A:glob-2} with $p=6$. Thus  in view of interpolation  
 we only have to obtain  
\begin{align}\label{4est}
\Big\| \prod_{k=1}^4
|A^\gamma[\psi_{\!J_k}]  {\mathcal P_{\!\mathbf {n}}}  P_{\!J_k}  f_k |^{\frac 14} 
\Big\|_{L^\frac83(\mathbb B^3(0,1) \times I )}
&\le C_{\varepsilon}
\lambda^{-\frac13+\varepsilon} \prod_{k=1}^4 
\|  f_{k} \|_{L^2(\mathbb R^3)}^{\frac 14}.
\end{align}

\begin{proof}[Proof of \eqref{4est}] For a given $\varepsilon>0$ we fix $\nu$ such that $10\nu=2^{-1}\vep$ and then take  an integer  $D$  such that 
$D\ge C_1/\nu$ with a large  constant $C_1$. 
For simplicity let us set  
\Be
\label{simple}
 F_k= A^\gamma[\psi_{\!J_k}]  {\mathcal P_{\!\mathbf {n}}}  P_{\!J_k}  f_k, \quad k=1, \dots, 4.
 \Ee
By Lemma \ref{op2} and \ref{op3}, we have 
\[
F_k
=F^{\bf c}_k + F^{+}_k+F^{-}_k + \mathcal E  f_{k}, \quad k=1, \dots, 4,
\]
where $\mathcal E$ satisfies $\| \mathcal E  f_k\|_q \le C_\delta \lambda^{ C-\nu D} \|f_k\|_p$ for $1\le p \le q \le \infty$, and 
\begin{align*}
F^{\bf c}_k &=\sum_{\substack{ |\ell | \le \lambda^{10\nu}}} e^{it\ell} 
\mathcal T^{\bf c}_\lambda (c_\ell \pi_{\bf c} {\widetilde \chi}_{\mathcal R_{J_k}} \widehat f_k), \\
F^{\pm}_k &=\sum_{0 \le m \le M-1} t^{-\frac {2m+1}2} \mathcal T^{\pm}_\lambda (\gamma_{m}^{\pm} \pi_{\bf o}^1  {\widetilde \chi}_{\mathcal R_{J_k}} \widehat f_k) .
\end{align*}

  We thus need to handle the product terms   $ \Pi_{k=1}^4 h_k$ where $h_k\in 
   \{ F^{\bf c}_k, F^{\pm}_k, \mathcal E  f_{k}\}$, $1\le k\le 4$. 
Any product which has $\mathcal E  f_{k}$ as one of its factors is easily handled by taking $C_1$ large enough  if one uses  H\"older's inequality and the trivial estimates    $\| \mathcal T_\lambda^{\bf c} (\pi_{\bf c}  \widehat g)\|_q\le C_\delta \lambda^{C}\|g\|_p$ and
$\| \mathcal T^{\pm}_\lambda  (\pi_{\bf o}^1 \widehat g)\|_q\le C_\delta \lambda^{C} \|g\|_p$, which hold for $1\le p\le q\le \infty $.   
So, it suffices to obtain the estimates for  the products which consist only  of   the terms $F^{\bf c}_k$, $F^{\pm}_k$. 
By \eqref{dec-m0}  and \eqref{dec-b1} we have  $\sum_{{ |\ell | \le \lambda^{10\nu}}}  \lambda^{\frac13-\nu}  \|c_\ell\|_\infty\le C\lambda^{3\nu} 
$ and $\sum_{\ell=0}^{M-1}  \|\gamma_{\ell}^{\pm}\|_\infty\le C_\delta \lambda^{-\frac13-\frac\nu 2}$. Thus, using the estimate   \eqref{T4} and Plancherel's theorem,  we 
obtain 
\[
\Big\| \prod_{k=1}^4
       |F^{\theta_k}_k |^{\frac 14} 
            \Big\|_{L^{\frac 83}(\mathbb B^4(0,2^3))}
                          \le C_\varepsilon \lambda^{-\frac13+10\nu+\frac \vep 2} \prod_{k=1}^4 \| f_k\|_2^{\frac 14},
\] 
where  $\theta_k\in\{\mathbf c, +,-\}$, $1\le k\le 4$.    Therefore we get \eqref{4est}. 
\end{proof}

\subsection{Proof of Proposition \ref{global}}
\label{sec:3.4}
We are in a position to prove Proposition \ref{global}.  
By Lemma \ref{localtoglobal}, it suffices   to show that 
\Be
\label{localmulti-2}
\Big\| \prod_{k=1}^4
|\widetilde \chi F_k |^{\frac 14} 
\Big\|_{L^p(\mathbb B^3(0,1)\times I,\omega)}
\le C_\delta \lambda^{-\varepsilon_p} 
\prod_{k=1}^4 \|  f_k\|_{L^p(\mathbb R^3)}^{\frac 14}
\Ee
for $14/5<p \le 6$. Here we keep using  the simpler notation \eqref{simple}.

We deduce the weighted estimate from Proposition \ref{prop:global} in the same way as in the proof of Proposition \ref{JC}. 
The difference is that we are dealing with a multilinear estimate and the exponent $p/4$ can be less than $1$.
Nonetheless, Lemma \ref{change} works as before. 
To  apply Lemma \ref{change},  we  break $\widetilde \chi F_k=\widetilde  A_kf_k + \mathcal E_kf_k$ where
\begin{align*}
\mathcal F( \widetilde  A_kf_k)(\xi,\tau)&= \beta_0 ((\lambda r_0)^{-1} \tau) \mathcal F(\widetilde \chi F_k)(\xi,\tau) 
\end{align*}
and $r_0=1+ 4\max\{|\gamma(s)|: s \in {\rm\supp} \psi_{J_k}, k=1,\dots, 4\} $.   
Since $[\omega]_3\le 1$ and $\|{\mathcal P_{\!\mathbf {n}}}  P_{\!J_k}  f\|_p\le C_\delta \|f\|_p$ and since $ | \mathcal E_k  f_k(x,t)| \le C\widetilde E_t^M \ast|{\mathcal P_{\!\mathbf {n}}}  P_{\!J_k}  f_k|(x)$ by Lemma \ref{lem:ker2}, we see  $ \| \mathcal E_k  f_k\|_{L^q(\mathbb R^3\times \mathbb R,\omega)} \le C_\delta \lambda^{-M} \|f_k\|_p$ for any $M>0$. Using the trivial  estimate $|\widetilde \chi F_k|\le  C_\delta \lambda^{3}(1+|\cdot|)^{-M}\ast| f_k|$, we also have 
$ \|\widetilde \chi F_k\|_{L^q(\mathbb R^3\times\mathbb R,\omega)} \le C_\delta \lambda^{3} \|f_k\|_p$. Making use of those estimates and taking a large $M$,  
one can easily see 
  \[
\Big\| \prod_{k=1}^4
|\widetilde \chi F_k |^{\frac 14} 
\Big\|_{L^q(\mathbb B^3(0,1)\times I,\omega)} \le C    
\Big\| \prod_{k=1}^4| \widetilde  A_k f_k |^{\frac 14} \Big\|_{L^q(\mathbb R^4,\omega)} + C_\delta \lambda^{-N} \prod_{k=1}^4 \|f_k \|_{L^p(\mathbb R^3)}^{\frac 14}\] 
for a large $N$ and $q\ge p$. 

 By  \eqref{control} we have $\| \prod_{k=1}^4
|\widetilde  A_k f_k |^{\frac 14} 
\|_{L^q(\mathbb R^4,\omega)}\le C  \lambda^{1/q} \| \prod_{k=1}^4| \widetilde  A_k f_k|^\frac14 \|_{L^q(\mathbb R^4)} $
since $[\omega]_3\le 1$ and the support of $\mathcal F(\prod_{k=1}^4 \widetilde  A_k f_k)$ is contained in a ball of radius $2^4r_0\lambda$. 
To estimate $\| \prod_{k=1}^4| \widetilde  A_k f_k|^\frac14 \|_{L^q(\mathbb R^4)}$, using the estimate  $ \| \mathcal E_k  f_k\|_{L^q(\mathbb R^3\times \mathbb R,\omega)} \le C_\delta \lambda^{-M} \|f_k\|_p$ again, we may disregard the minor contributions. So, it is sufficient to  consider  $ \|  \prod_{k=1}^4 
|\widetilde \chi F_k |^{\frac 14}\|_{L^q(\mathbb R^4)}.$
Since $\supp \widetilde \chi \subset I$,  by the estimate \eqref{A:glob} we get
\[
\Big\| \prod_{k=1}^4
|\widetilde \chi F_k |^{\frac 14} 
\Big\|_{L^q(\mathbb B^3(0,1)\times I,\omega)}	\le
C_{\vep}(\delta)
\lambda^{\frac 7{24}(\frac1{p}-\frac 5{14})+\vep}
\prod_{k=1}^4 \|f_k \|_{L^p(\mathbb R^3)}^\frac14
\]
 for 
 $1/q=5/(8p)+1/16$
 and $2\le p \le 6$. Finally, we obtain \eqref{localmulti-2}  for $14/5<p \le 6$ by  H\"older's inequality since $q\ge p$ and $\|\omega\|_{L^1(\mathbb B^3(0,1)\times I)}\le C [\omega]_3$. 
\qed

\begin{rmk} 
\label{order-} In the above  we try to obtain the estimate \eqref{gb} on a range of $p$ as large as possible by suppressing  $\nu$ arbitrarily  small (\emph{Proof of \eqref{4est}}).  This forces us to take a large  $D\ge C_1/\nu$. However, to obtain the maximal estimate   it is enough to have the estimate \eqref{gb} on a smaller range $3<p\le 6$  instead of $14/5<p \le 6$.  For $3<p\le 6$, we can prove \eqref{gb} with a fixed $\nu$ and $D$. For example, optimizing the estimates at various places, we can take
$\nu=1/397$ and $D=720$. In other words, Theorem \ref{max} holds true for $\gamma\in \mathrm C^{720}(\mathbb J)$.
\end{rmk}

\section{Proof of Theorem \ref{max}}\label{sec4}
In this section we complete the proof of Theorem \ref{max}. We prove the sufficiency  and the necessity parts in separate sections.  

\subsection{Sufficiency} 
\label{s-conclusion}
By the reduction in Section \ref{sec2.4}, Lemma  \ref{WtoM} and Lemma \ref{measure-max}, 
it suffices to prove Proposition \ref{QQ}, which also proves Theorem \ref{LSA}.  

\subsubsection*{Decomposition}  We first  decompose the averaging operator $A^\gamma[\psi]$ in such a  way that we can use the multilinear estimate obtained in Section \ref{sec3}.     The following Lemma 
\ref{decomp-lem} is a slight modification of \cite[Lemma 2.8]{HL}.  Let us set
\begin{align*}
\mathfrak J^4_\ast(\delta)=\big\{ (J_1,\dots, J_4) : J_1, \dots, J_4 \in \mathfrak J_\circ(\delta),
\quad~\quad~ \min_{\ell \neq k} \dist (J_\ell, J_k) \ge   \delta \big\}.
\end{align*}

\begin{lem}\label{decomp-lem} Let  $\psi\in \mathfrak N^D(J_\circ)$ and $\gamma \in \mathfrak C^D(\varepsilon_\circ)$. 
There is a constant  $C=C(D)$ independent of $z=(x,t)$, $\gamma$, and $\delta$  such that
\begin{equation}
\label{4d}
| A^\gamma[\psi] f(z)|
\le
C\!\!\! \max_{J \in \mathfrak J_\circ(\delta)}| A^\gamma[\psi_{\!J}]   f(z)|   
+C\delta^{-1} \!\!\!\!\!\!\sum_{(J_1,\dots,J_4) \in \mathfrak J^4_\ast(\delta)}
 \prod_{k=1}^4  \big|A^\gamma[\psi_{\!J_k}]  f(z) \big|^{\frac 14},
\end{equation}
where $\psi_{\!J}\in \mathfrak  N^D(J)$ and $\psi_{\!J_k}\in \mathfrak  N^D(J_k)$. 
\end{lem}

\begin{proof} Let us recall  \eqref{Jsum}. It is clear that  there is a constant $C_D>0$ such that 
$C_D^{-1}\psi \zeta_J \in \mathfrak  N^D(J)$ for $J\in \mathfrak J_\circ(\delta)$. 
Setting  $\psi_{\!J}  =C_D^{-1}\psi \zeta_J$ we  have   
\[A^\gamma [\psi] f(z)= C_D \sum_{J\in \mathfrak J_\circ(\delta)} A^\gamma [\psi_{\!J}]  f(z).\] 
 
 Let us set  $\mathfrak J_1=\mathfrak J_\circ(\delta)$. For a fixed $z$, define $J_1^*$ to be an interval  in $\mathfrak J_1$ such that 
$
|A^\gamma[\psi_{J_1^*}]  f(z)|=\max_{J \in \mathfrak J_1} | A^\gamma[\psi_{\!J}]  f(z)|.
$
For $k=2,3,4,$  we recursively define $\mathfrak J_k$ and $J_k^*\in \mathfrak J_k$.  Let  $
\mathfrak J_k= \{J \in \mathfrak J_{k-1}: \dist(J, J_{k-1}^*) \ge \delta \}
$
and let $J_k^*\in \mathfrak J_k$ denote an interval such that $|A^\gamma[\psi_{\!J_k^*}] f(z)|= \max_{J \in \mathfrak J_k } |A^\gamma[\psi_{\!J}]  f(z)|.$
Thus, if
$\dist(J,J_{k}^*)\ge \delta$ for all $1\le k\le 3$, we have $|A^\gamma[\psi_{\!J}]   f|\le |A^\gamma[\psi_{\!J_k^*}]  f|$ for $1\le k\le 4$.

Let us denote $\mathcal J=\bigcup_{k=1}^3\{J\in \mathfrak J_\circ(\delta): \dist(J, J_k^*) < \delta \}$. 
Splitting the sum into the cases  $J\in \mathcal J$ and  $J\not\in  \mathcal J$, we have
\begin{align*}
C_D^{-1}|A^\gamma[\psi] f(z) | \le 
\sum_{J\in \mathcal J} 
|A^\gamma[\psi_{\!J}]   f(z) |
+
\sum_{J\not\in \mathcal J} |A^\gamma[\psi_{\!J}]   f(z) |  .
\end{align*}
The first sum on the right hand side is apparently bounded by a constant times $\max_{J \in \mathfrak J_\circ(\delta)} |A^\gamma[\psi_{\!J}]   f(z)|$ and the second  by $ C\delta^{-1} \prod_{k=1}^4 | A^\gamma[\psi_{\!J_k^*}]  f(z)|^{\frac 14}.$ 
This gives \eqref{4d} since $\dist( J_k^*, J^*_\ell)\ge \delta$ if $k\neq \ell$.
\end{proof}

In the next lemma, using $K_\delta(\lambda)$ given in Lemma \ref{scaling2}, we get a bound on the first one on the right hand side of \eqref{4d}. 

\begin{lem}\label{lem:lin}  Let $2< p \le 6$, and let  $[\omega]_3\le 1$ and   $\psi_{\!J}\in \mathfrak N^D(J)$   for each  $J\in \mathfrak J_\circ(\delta)$.
If $\delta^3\lambda\ge 2^2$ and $\varepsilon_\circ>0$ is sufficiently small,  there 
is an  $\vep_p>0$ such that  
\begin{align*}
\|\max_{J \in \mathfrak J_\circ(\delta)} |A^\gamma[\psi_{\!J}]   f|\big\|_{L^p(\mathbb R^3 \times [1,2],\omega)} 
&\le 
C
\big( \delta^{1-\frac {3}p} K_\delta(\lambda)
+  C_\delta \lambda^{-\vep_p} \big)\|f\|_{L^p(\mathbb R^3)}
\end{align*}
holds whenever  $\gamma \in \mathfrak C^D(\varepsilon_\circ)$ and $\widehat f$ is supported on $\mathbb A_\lambda$.  
\end{lem}

\begin{proof}[Proof of Lemma \ref{lem:lin}]
By the embedding $\ell^p \subset \ell^\infty$
and Minkowski's inequality,
\[
\|  \max_{J \in \mathfrak J_\circ(\delta)} 
|A^\gamma {[\psi_{\!J}]} f|  \big\|^p_{L^p(\mathbb R^3 \times [1,2],\omega)}
\le 2^p( \mathrm{I}+  \mathrm{I\!I}),\]
where
\begin{align*}
\mathrm{I}&= \sum_{J \in \mathfrak J_\circ(\delta)} 
\| A^\gamma[\psi_{\!J}] P_{\!J} f \big\|_{L^p(\mathbb R^3\times [1,2],\omega)}^p, 
\\
 \mathrm{I\!I}&=\sum_{J \in \mathfrak J_\circ(\delta)} 
\| A^\gamma[\psi_{\!J}]  (f-P_{\!J} f) \big\|_{L^p(\mathbb R^3\times  [1,2],\omega)}^p.
\end{align*}

For  
$\rm{I\!I}$ we apply   Proposition  \ref{JC}.  Taking $\vep_p=\frac14(\frac12-\frac1p)$ and 
using the estimate \eqref{small} with $\vep=\varepsilon_p/2$, we have
$
\mathrm{I\!I}^{\,\frac1p} \le   C_\delta 
\lambda^{-\varepsilon_p}  \|f\|_{L^p(\mathbb R^3)} 
$ 
since there are at most $C\delta^{-1}$ many $J$.
To handle $\rm I$, we invoke Lemma \ref{scaling2} and  then use  Lemma \ref{ortho} to obtain
\[\mathrm I 
\le C 
\delta^{p-3}  
K_\delta(\lambda)^p 
\sum_{J \in \mathfrak J_\circ(\delta)} \|P_{\!J}f\|_p^p
\le  C 
\delta^{p-3}  
K_\delta(\lambda)^p 
 \|f\|_p^p.\]
   Therefore  the desired bound  follows.
\end{proof}

Now we consider  the  product  terms appearing in \eqref{4d}.

\begin{lem}\label{lem:multi} Let $\frac{14}5 <p \le 6$, $[\omega]_3\le 1$, and  $(J_1,\dots,J_4)\in \mathfrak J^4_\ast(\delta)$.
If  $\delta^3\lambda\ge 2^2$  and $\varepsilon_\circ>0$ is small enough,
there are positive constants $\varepsilon_p$, $c$, $D$  such that 
\begin{align}\label{eq:multi}
\Big\|  
\prod_{k=1}^4 |A^\gamma[\psi_{\!J_k}] f|^{\frac 14} \Big\|_{L^p(\mathbb R^3\times  [1,2],\omega)} \!\!\le
C_\delta     \big(\lambda^{-\vep_p}+  \lambda^{-c} 
K_\delta(\lambda) \big) \|f\|_{L^p(\mathbb R^3)}
\end{align}
holds whenever $\gamma \in \mathfrak C^D(\varepsilon_\circ)$, $\psi_{J_k}\in \mathfrak N^D(J_k)$, $k=1,\dots, 4,$ and $\widehat f$ is supported in $\mathbb A_\lambda$.
\end{lem}

\begin{proof}
For each $1\le k\le 4$  we split   
$f= b_k+g_k$, 
 where 
 \[ b_k= {\mathcal P_{\!\mathbf {n}}} P_{J_k}f, \quad g_k= {\mathcal P_{\!\mathbf {n}}}(1-P_{\!J_k}) f+ {\mathcal P_{\!\mathbf {e}}} f.\] 
We here use $f={\mathcal P_{\!\mathbf {n}}} f+{\mathcal P_{\!\mathbf {e}}} f$ because  $\widehat f$ is supported on $\mathbb A_\lambda$.  Thus, 
the left hand side of \eqref{eq:multi} is bounded by 
a constant times 
\begin{align*}
\mathfrak M=\sum_{h_k\in \{b_k, g_k\} }
\Big\| 
\prod_{k=1}^4 
|A^\gamma[\psi_{\!J_k}] h_{k} |^{\frac 14} \Big\|_{L^p(\mathbb R^3 \times  [1,2],\omega)}.
\end{align*}
We consider the  cases $(h_1,\dots, h_4)=(b_1,\dots, b_4)$ and $(h_1,\dots, h_4)\neq (b_1,\dots, b_4)$.  For  the former case we use  Proposition \ref{global}  and the estimate \eqref{easy2}. Since $14/5<p\le 6$, there is an $\vep_p>0$ such that 
\begin{align*}
\Big\|\prod_{k=1}^4 
|A^\gamma[\psi_{\!J_k}] b_{k} |^{\frac 14} \Big\|_{L^p(\mathbb R^3 \times  [1,2],\omega)}
\le C_\delta 
\lambda^{-\varepsilon_p}
\|f\|_{L^p(\mathbb R^3)}. 
\end{align*}
For the other case we combine Proposition \ref{JC}, \ref{JC*},  and Lemma \ref{scaling2}. In fact, 
Proposition \ref{JC} and \ref{JC*}  followed by \eqref{easy2}
yield
\[
\| A^\gamma[\psi_{\!J_k}] g_{k} \|_{L^p(\mathbb R^3 \times [1,2], \omega)}
\le C_\varepsilon \delta^{-C} \lambda^{\frac 12(\frac1p-\frac12)+\varepsilon}\|f\|_{L^p(\mathbb R^3)}
\]
 for $2\le p\le 6$.  If we consider a particular case $(h_1,\dots, h_4)= (b_1,b_2, b_3, g_4)$, by  H\"older's inequality and the above estimate    we have 
\[
\Big\| 
\prod_{k=1}^4 
|A^\gamma[\psi_{\!J_k}] h_{k} |^{\frac 14} \Big\|_{L^p(\mathbb R^3 \times  [1,2],\omega)}\!\!\le  C_\delta\lambda^{-c}\|f\|_{L^p(\mathbb R^3)}^\frac14\!\!\prod_{k=1}^3
\big\|  
A^\gamma[\psi_{\!J_k}] b_{k} \big\|_{L^p(\mathbb R^3\times  [1,2],\omega)}^{\frac 1{4}}\] 
for a constant $c>0$ because $p> 14/5$.  We apply  Lemma \ref{scaling2} to handle the  last three factors. Since $\|  b_k \|_{L^p(\mathbb R^3)} \le C_1 \delta^{-C}  \|f\|_{L^p(\mathbb R^3)}$ from \eqref{easy2},   the inequality \eqref{hoho}  gives 
\[
\Big\| 
\prod_{k=1}^4 
|A^\gamma[\psi_{\!J_k}] h_{k} |^{\frac 14} \Big\|_{L^p(\mathbb R^3 \times  [1,2],\omega)}\le  C_\delta  \lambda^{-c}
K_\delta(\lambda)^\frac34  \|f\|_{L^p(\mathbb R^3)}.\]
We can deal with the remaining products similarly.  As a consequence, we obtain 
\[\mathfrak M\le C_\delta    \Big(\lambda^{-\vep_p}+  \sum_{\ell=1}^3 \lambda^{-(4-\ell) c} 
K_\delta(\lambda)^\frac\ell4 \Big) \|f\|_{L^p(\mathbb R^3)} \] 
and therefore the bound \eqref{eq:multi} after a simple manipulation since we may assume $\vep_p\le c$ taking a smaller $\vep_p$ if necessary. \end{proof}

 We  now conclude the proof of \eqref{Q} putting together the  previous estimates. 

\subsubsection*{Proof of  \eqref{Q}}  Since   $\|Af\|_{L^\infty(\mathbb R^3 \times [1,2], \omega)}\le C\|f \|_{L^\infty(\mathbb R^3)}$,   by interpolation it is sufficient to show \eqref{Q} for $3<p<6$. Let $p\in (3,6)$ and take an $\vepc>0$  small  enough and a large $D$ such that  the estimates in  Lemma \ref{lem:lin} and \ref{lem:multi}  hold whenever 
$\gamma \in \mathfrak C^D(\varepsilon_\circ)$ and $\psi_{\!J}\in \mathfrak N^D(J)$, $J\in \mathfrak J_\circ (\delta)$.

Let $\gamma \in \mathfrak C^D(\varepsilon_\circ),$ $\omega\in \Omega^3 $ with $[\omega]_3\le 1$ and $\psi\in  \mathfrak N^D(J_\circ)$, and let $f$ be a function such that  $\supp \widehat f \subset \mathbb A_\lambda$ and $\|f\|_{p}\le 1$.  By \eqref{4d}  and Minkowski's inequality  we see that 
$\| A^\gamma[\psi]   f\|_{L^p(\mathbb R^3\times [1,2],\omega)}$ is bounded by 
\begin{align*}
C\big\| \max_{J \in \mathfrak J_\circ(\delta)} 
|A^\gamma[\psi_{\!J}]  f |\big\|_{L^p(\mathbb R^3\times [1,2],\omega)} +C_\delta \!\!\!\! \!\!\sum_{(J_1,\dots,J_4) \in \mathfrak J^4_\ast(\delta)}
\Big\| 
\prod_{k=1}^4 |A^\gamma[\psi_{\!J_k}]f|^{\frac 14} \Big\|_{L^p(\mathbb R^3\times [1,2],\omega)}. 
\end{align*}
Then Lemma \ref{lem:lin} and \ref{lem:multi} give
\[ \| A^\gamma [\psi]   f\|_{L^p(\mathbb R^3\times [1,2],\omega)} \le C
\big( \delta^{1-\frac {3}p} +\lambda^{-c} \big) K_\delta(\lambda)
+ C_\delta \lambda^{-\vep_p}\] 
 if $2^{2}\delta^{-3}\le \lambda$. Taking supremum over $f, \omega, \psi$, and $\gamma$,
we obtain 
\begin{equation}\label{eq11}
\begin{aligned}
Q(\lambda)
\le 
           C \big( \delta^{1-\frac {3}p} +\lambda^{-c} \big) K_\delta(\lambda)
                    + C_\delta \lambda^{-\vep_p}
\end{aligned}
\end{equation}
for $2^{2}\delta^{-3}\le \lambda.$
In order to close the induction we  need to modify $Q(\lambda)$ slightly. Fix $0<b$, which is to be chosen later. We define 
\begin{align*}
\overline Q_b(\lambda)= \sup_{1 \le r \le \lambda}
r^{b} Q(r).
\end{align*}

We observe  $\lambda^b K_\delta(\lambda)
   \le  2^{2b} \delta^{-3b} 
\sum_{2^{-2}\delta^3 \lambda \le 2^{j} \le 2^2\delta\lambda } 2^{jb} Q(2^{j})$, and  hence  we have 
$
\lambda^b K_\delta(\lambda)
   \le C |\log\delta| \delta^{-3b} \overline Q_b(2^2\delta\lambda).
$
Multiplying $\lambda^b$ to both sides of  \eqref{eq11},  we get  
\[ \lambda^b Q(\lambda)\le C
\big( \delta^{1-\frac {3}p} +\lambda^{-c} \big) |\log\delta| \delta^{-3b} \overline Q_b(2^2\delta\lambda)+ C_\delta \lambda^{b-\vep_p}
\] 
for $2^{2}\delta^{-3}\le \lambda$.
We now choose a small  $b$ such that $1-\frac {3}p-3b>0$ and $b-\vep_p<0$, then fix  a small enough $\delta>0$   such that $C \delta^{1-\frac {3}p}  |\log\delta| \delta^{-3b}\le 2^{-2}$ and  $2^2\delta\le 1$. Such a choice is clearly possible because $p>3$.  Let $\lambda_\circ$ be a large number such that $\delta^{1-\frac {3}p}\ge \lambda^{-c}_\circ$ and $2^{2}\delta^{-3}\le \lambda_\circ$. Then 
we have the inequality  $ \lambda^b Q(\lambda)\le 2^{-1} \overline Q_b(\lambda)+ C_\delta$  for $\lambda \ge  \lambda_\circ$ since $\overline Q_b$ is increasing. 
This obviously implies 
\[ \lambda^b Q(\lambda)\le 2^{-1} \overline Q_b(r)+ C_\delta
\] 
for $\lambda_\circ\le \lambda\le r$.  Note $\overline Q_b(\lambda_\circ)\le \lambda_\circ^b C_2$ for some constant $C_2$ (because of 
the trivial estimate $Q(\lambda)\le C\lambda^2$). Taking supremum over $\lambda\in [1, r]$  we get   $\overline Q_b(r)\le 2^{-1} \overline Q_b(r)+ \lambda_\circ^b C_2+ C_\delta$. Therefore we have  $\overline Q_b(r)\le C_3$ for a constant $C_3$ and  conclude   $Q(\lambda)\le C_3\lambda^{-b}$ for $\lambda\ge 1$.   \qed

\begin{rmk} 
\label{further-result} Routine adaptation of our argument also proves $L^p$ improving property of the localized maximal operator $\overline Mf(x):= \sup_{1\le t\le 2} \big| A  f(x,t) \big|$. In fact,  the estimate $
\|\overline Mf\|_{L^q(\mathbb R^3)}
\le C \| f\|_{L^p(\mathbb R^3)}$ holds provided that $(1/p,1/q)$ is contained in the interior of the triangle with vertices $(0,0), (1/3,1/3),$ and $(19/66, 8/33)$.  
It is possible to  extend the range slightly making use of the estimate  \eqref{ps14*} for $p>6$. Furthermore, one can show that 
  $\overline M$ is bounded from $L^p$ to $L^p(d\mu)$  for $p>9-2\alpha$ when $\mu$ is an  $\alpha$ dimensional measure  and  $3> \alpha>\frac{65- \sqrt{865}}{12}=2.9657\dots.$  
  \end{rmk}

\subsection{Necessity}\label{Sec:N} 
To prove that $L^p$ boundedness of $M$ fails for $p\le 3$, 
it is sufficient to show the next proposition.   
Our construction below is a  modification of  Stein's example in \cite{Stein2}.

\begin{prop} Let $p\le 3$ and  $\psi \not\equiv 0$ be a nonnegative continuous function supported in $\mathbb J$. Suppose $\gamma : \mathbb J\rightarrow \mathbb R^3$ is a smooth curve  
satisfying \eqref{nonv}.  Then there is an $h\in L^p(\mathbb R^3)$  such that $Mh =\infty$ on a nonempty open set.
\end{prop}

\begin{proof}
Since $\psi \ge0$ and $\psi \not\equiv0$,  we may assume that $\psi(s)\ge c$  on an interval $J \subset \mathbb J$ for some $c>0$. By \eqref{nonv} 
we may additionally assume that $|\gamma(s)|\ge c$
on  $J$  taking a subinterval of $J$  if necessary  because the condition  \eqref{nonv} can not be satisfied if there is no such a subinterval.

Since $\gamma'(s),$ $ \gamma''(s),$ $ \gamma'''(s)$ are linearly independent, 
we can write 
\begin{align}\label{GT}
\gamma(s) 
=c_1(s)\gamma'(s)+c_2(s)\gamma''(s)+c_3(s)\gamma'''(s), \quad s\in J
\end{align}
for some  smooth functions $c_1$, $c_2,$ and $c_3$.
We claim that there is an $s_\circ \in J$ such that $c_3(s_\circ)\neq0$.
Suppose that there is no such  $s_\circ \in J$,  that is to say,  $c_3(s) \equiv 0$ for all $s \in J$.
Differentiating both  side of  \eqref{GT},  we have
$
(c_1'(s)-1)\gamma'(s)+[c_1(s)+c_2'(s)]\gamma''(s)+c_2(s)\gamma'''(s)=0,
$
which  implies $c_2(s) \equiv 0$,    $c_1(s)+c_2'(s) \equiv 0$, and   $c_1'(s) \equiv 1$  for $s \in J$. This leads to a contradiction and proves the claim. 
Therefore there are  $s_\circ\in J$ and $\delta>0$ such that 
\[
|c_3(s)| \ge c, \quad s\in [s_\circ-\delta, s_\circ+\delta] \subset J
\]
for some   $c>0$.  We only consider the case $c_3(s)\ge c$ since the other case can be handled similarly.

For $x \in \mathbb R^3$  let $y=(y_1, y_2, y_3)$  denote the coordinate of $x$ with respect to the basis  $\{ \gamma'(s_\circ), \gamma''(s_\circ), \gamma'''(s_\circ)\}$, i.e., 
$
x=y_1 \gamma'(s_\circ)+y_2 \gamma''(s_\circ)+y_3\gamma'''(s_\circ), 
$
and set $\overline{y}=(y_1, y_2)$. For some $\varepsilon\in (0,1/3)$ we take  $g(t)= \chi_{[0,2^{-1}]}(t) |t|^{-\frac 13} 
|\log |t||^{-\frac 13-\varepsilon}$ and then we consider 
\[
h(x)=\chi_0(|\overline {y}|) g(y_3),
\]
where $\chi_0 \in \mathrm C_0^\infty ([-2,2])$ is a nonnegative function such that  $\chi_0=1$ on
$[-1,1]$. It is easy to see $h\in L^p(\mathbb R^3)$ for $p\le 3$ because  $g\in  L^p(\mathbb R)$ for $p\le 3$.  
Thus we only have to show  that  $\sup_{0<t}A  h=\infty$ on a nonempty  open set.

We  write  
$
\gamma(s)=  a_1(s) \gamma'(s_\circ)+ a_2(s)\gamma''(s_\circ) +a_3(s) \gamma'''(s_\circ)
$ 
and $\overline a(s)=(a_1(s),a_2(s))$.   Since $c_j(s_\circ)=a_j(s_\circ)$, $j=1,2,3$,  by a Taylor  expansion   we have
\begin{equation*}
\begin{aligned}
\gamma(s)
&=
\big( c_1(s_\circ)+(s-s_\circ) \big) \gamma'(s_\circ)
+\big( c_2(s_\circ)+ (s-s_\circ)^2/2! \big) \gamma''(s_\circ)
\\
&\qquad \ +\big( c_3(s_\circ)+(s-s_\circ)^3/3! \big) \gamma'''(s_\circ)
+\mathcal O\big((s-s_\circ)^4\big).
\end{aligned}
\end{equation*}
So,  $y_3-ta_3(s)=y_3-t c_3(s_{\circ})-t\big( 6^{-1}(s-s_\circ)^3+\mathcal O((s-s_\circ)^4) \big)$. For $y_3>0$ we  take   $t=y_3/c_3(s_\circ)>0$. 
Then  it follows that $C_1y_3 |s-s_\circ|^3\le  |y_3-ta_3(s)|\le C_2
 y_3 |s-s_\circ|^3
$ 
for some $C_1, C_2>0$, so   $|g(y_3-ta_3(s))|\ge C y_3^{-\frac 13} |s-s_\circ|^{-1}
| \log (y_3 |s-s_\circ|^3) |^{-\frac 13-\varepsilon}$ 
 provided that  $|s-s_\circ|<c'$ for a small $c'>0$ and $0<y_3\le 1$. 
 Thus,  by our choice of $\delta$ and $s_\circ$   we have 
\begin{align*}
A  h\Big(x,\frac{y_3}{c_3(s_\circ)}\Big)
&\ge C 
y_3^{-\frac 13}
\int_{|s-s_\circ|\le  \delta' }  
\widetilde 
\chi_0 (y,s) |s-s_\circ|^{-1}
| \log (y_3 |s-s_\circ|^3) |^{-\frac 13-\varepsilon} \,ds
\end{align*} 
for  $0<y_3\le 1$  where  $\delta'=\min(\delta,c')$ and  $\widetilde 
\chi_0 (y,s)=\chi_0 (| \overline{y}
-\frac{y_3}{c_3(s_\circ)}  \overline a(s) |)$.   Since  $\widetilde 
\chi_0 (y,s)\ge 1$ if $|y|\le r_\circ$ for a small enough $r_\circ>0$, 
we have  \begin{align*}
A  h\Big(x,\frac{y_3}{c_3(s_\circ)}\Big)
&\ge C
y_3^{-\frac 13}\! \int_{|s-s_\circ| \le  \min(\delta',\,  y_3^{\frac 13}/10)}
|s-s_\circ|^{-1}
| \log |s-s_\circ| |^{-\frac 13-\varepsilon} \,ds =\infty
\end{align*}
for $y\in \mathbb B^3(0, r_\circ)\cap \{y: 0<y_3<1\}$  
as desired. 
\end{proof}

\subsection*{Acknowledgement}  
This work was supported by  the National Research Foundation of Korea  (NRF)  grants  number NRF-2019R1A6A3A01092525 (Hyerim Ko) and  NRF-2021R1A2B5B02001786 (Sanghyuk  Lee and Sewook  Oh).

\bibliographystyle{plain}

\end{document}